\documentclass[12pt,reqno,dvipdfmx]{amsart}

\usepackage{amsmath,amssymb,amsthm,mathtools}
\mathtoolsset{showonlyrefs=true}
\usepackage{mathrsfs}
\usepackage{bm}
\usepackage{enumerate}

\numberwithin{equation}{section}
\allowdisplaybreaks[2]

\usepackage[dvipdfmx]{xcolor}

\usepackage{aliascnt}

\theoremstyle{plain}

\newtheorem{theorem}{Theorem}[section]

\newaliascnt{proposition}{theorem}
\newtheorem{proposition}[proposition]{Proposition}
\aliascntresetthe{proposition}

\newaliascnt{lemma}{theorem}
\newtheorem{lemma}[lemma]{Lemma}
\aliascntresetthe{lemma}

\newaliascnt{corollary}{theorem}
\newtheorem{corollary}[corollary]{Corollary}
\aliascntresetthe{corollary}

\newaliascnt{conjecture}{theorem}

\aliascntresetthe{conjecture}

\theoremstyle{definition}

\newaliascnt{definition}{theorem}
\newtheorem{definition}[definition]{Definition}
\aliascntresetthe{definition}

\newaliascnt{example}{theorem}
\newtheorem{example}[example]{Example}
\aliascntresetthe{example}

\newaliascnt{notation}{theorem}

\aliascntresetthe{notation}

\theoremstyle{remark}

\newaliascnt{remark}{theorem}

\aliascntresetthe{remark}

\newaliascnt{question}{theorem}
\newtheorem{question}[question]{Question}
\aliascntresetthe{question}

\usepackage[
  dvipdfmx,
  unicode=true,
  colorlinks=true,
  linkcolor=blue,
  citecolor=blue,
  urlcolor=blue
]{hyperref}

\usepackage[nameinlink,noabbrev]{cleveref}

\crefname{theorem}{Theorem}{Theorems}
\Crefname{theorem}{Theorem}{Theorems}

\crefname{proposition}{Proposition}{Propositions}
\Crefname{proposition}{Proposition}{Propositions}

\crefname{lemma}{Lemma}{Lemmas}
\Crefname{lemma}{Lemma}{Lemmas}

\crefname{corollary}{Corollary}{Corollaries}
\Crefname{corollary}{Corollary}{Corollaries}

\crefname{conjecture}{Conjecture}{Conjectures}
\Crefname{conjecture}{Conjecture}{Conjectures}

\crefname{definition}{Definition}{Definitions}
\Crefname{definition}{Definition}{Definitions}

\crefname{example}{Example}{Examples}
\Crefname{example}{Example}{Examples}

\crefname{notation}{Notation}{Notations}
\Crefname{notation}{Notation}{Notations}

\crefname{remark}{Remark}{Remarks}
\Crefname{remark}{Remark}{Remarks}

\crefname{question}{Question}{Questions}
\Crefname{question}{Question}{Questions}

\crefname{section}{\S}{\S\S}
\Crefname{section}{\S}{\S\S}

\crefname{subsection}{\S}{\S\S}
\Crefname{subsection}{\S}{\S\S}

\crefname{subsubsection}{\S}{\S\S}
\Crefname{subsubsection}{\S}{\S\S}

\crefname{appendix}{Appendix}{Appendices}
\Crefname{appendix}{Appendix}{Appendices}

\newcommand{\Z}{\mathbb Z}
\newcommand{\N}{\mathbb N}
\newcommand{\R}{\mathbb R}
\newcommand{\C}{\mathbb C}
\newcommand{\Q}{\mathbb Q}

\newcommand{\cN}{\mathcal N}

\newcommand{\fkm}{\mathfrak m}

\newcommand{\fkJ}{\mathfrak J}

\newcommand{\fkM}{\mathfrak M}

\newcommand{\fkP}{\mathfrak P}
\newcommand{\fkQ}{\mathfrak Q}

\newcommand{\fkT}{\mathfrak T}
\newcommand{\fkU}{\mathfrak U}

\newcommand{\ba}{{\boldsymbol a}}
\newcommand{\bb}{{\boldsymbol b}}

\newcommand{\bs}{{\boldsymbol s}}
\newcommand{\bt}{{\boldsymbol t}}
\newcommand{\bu}{{\boldsymbol u}}
\newcommand{\bv}{{\boldsymbol v}}
\newcommand{\bw}{{\boldsymbol w}}
\newcommand{\bx}{{\boldsymbol x}}
\newcommand{\by}{{\boldsymbol y}}
\newcommand{\bz}{{\boldsymbol z}}
\newcommand{\bzero}{{\boldsymbol 0}}

\newcommand{\aalpha}{{\boldsymbol{\alpha}}}
\newcommand{\bbeta}{{\boldsymbol{\beta}}}
\newcommand{\ggamma}{{\boldsymbol{\gamma}}}

\newcommand{\nnu}{{\boldsymbol{\nu}}}
\newcommand{\bttheta}{{\boldsymbol{\theta}}}
\newcommand{\xxi}{{\boldsymbol{\xi}}}

\DeclareMathOperator{\supp}{supp}

\DeclareMathOperator{\ord}{ord}

\DeclareMathOperator{\Sol}{Sol}

\DeclareMathOperator{\Ann}{Ann}

\DeclareMathOperator{\Ker}{Ker}

\DeclareMathOperator{\ini}{in}
\DeclareMathOperator{\ind}{ind}
\DeclareMathOperator{\find}{find}

\DeclareMathOperator{\nsupp}{nsupp}
\DeclareMathOperator{\NS}{NS}

\title{Logarithmic $A$-hypergeometric series 
${\textrm I}\! {\textrm I}\! {\textrm I}$}
\keywords{$A$-hypergeometric systems; hypergeometric series; Frobenius's method.}
\subjclass[2020]{Primary: 33C70}

\author{Go Okuyama}
\address{
\begin{flushleft}
(Go Okuyama)\\
Higher Education Support Center\\
Hokkaido University of Science\\
Sapporo, 006-8585, Japan.
\end{flushleft}
}
\email{gokuyama@hus.ac.jp}

\author{Mutsumi Saito}
\address{
\begin{flushleft}
(Mutsumi Saito)\\
Department of Mathematics, Faculty of Science\\
Hokkaido University\\
Sapporo, 060-0810, Japan.
\end{flushleft}
}
\email{saito@math.sci.hokudai.ac.jp}

\begin{document}
\begin{abstract}
We study the logarithmic coefficients that can occur at a fixed fake exponent in an $A$-hypergeometric series subject to prescribed negative support conditions.

Let $L=\Ker_{\mathbb Z}(A)$.
For a fixed generic weight $\bw$, a fixed fake exponent $\bv_0$,  and an ordered negative support family, we first derive a finite system of constant-coefficient differential equations whose solutions encode the admissible logarithmic coefficients.
A normalization of the coefficient equations shows that, for $\bu\in L$, the normalized coefficient associated with $\bx^{\bv_0+\bu}$ depends only on the negative support of $\bv_0+\bu$.
From this finite system, we identify the annihilator of the coefficient space with an explicitly defined colon ideal.
For the negative support family determined by the direction $\bw$, we further identify this colon ideal with the primary component of the indicial ideal supported at $\bv_0$, shifted to the origin.

We next introduce an ambient perturbation construction in which the fake exponent is perturbed in the full ambient space rather than only within the affine space $\bv_0+L_{\mathbb C}$.
We prove that the ambient perturbation construction produces $A$-hypergeometric series and realizes the full coefficient space.

Finally, we compare the ambient perturbation construction with the intrinsic perturbation construction developed in our previous papers.
The intrinsic construction always yields a subspace of the full coefficient space, and it realizes the full coefficient space if and only if a natural equality between the corresponding colon ideals holds.
We also give several sufficient conditions for this equality.
\end{abstract}
\maketitle

\section{Introduction}

This is the third paper in a series \cite{Log,Log2} devoted to Frobenius's method for $A$-hypergeometric systems.
The purpose of the present paper is to determine, under prescribed negative support conditions, the logarithmic coefficients that can occur at a fixed fake exponent in an $A$-hypergeometric series.
The problem is formulated for an ordered family of negative supports and is reduced to a finite system of constant-coefficient differential equations.
The resulting coefficient space is then described as the orthogonal complement of an explicit colon ideal.
We also give an ambient perturbation construction that realizes the full coefficient space and compare it with the intrinsic perturbation construction developed in \cite{Log,Log2}.

We begin by recalling the basic notation.
Let $A=(\ba_1,\dots,\ba_n)=(a_{ij})$ be a $d\times n$ integer matrix of rank $d$.
Throughout this paper we assume that $A$ is homogeneous, meaning that all columns $\ba_j$ lie in a single affine hyperplane of $\Q^d$ not passing through the origin.
Set
\[
L:=\Ker_{\mathbb Z}(A)\subset\mathbb Z^n,\quad L_{\mathbb C}:=\mathbb C\otimes_{\mathbb Z}L=\Ker_{\mathbb C}(A)\subset\mathbb C^n.
\]
Let $\N$ denote the set of nonnegative integers.
The toric ideal associated with $A$ is
\begin{equation}
\label{eq:ToricIdeal}
I_A=\left\langle\partial_{\bx}^{\bu}-\partial_{\bx}^{\bv}\,\middle|\,
A\bu=A\bv,\ \bu,\bv\in\N^n\right\rangle\subset\C[\partial_{\bx}],
\end{equation}
where $\C[\partial_{\bx}]=\C[\partial_1,\dots,\partial_n]$.
For a parameter $\bbeta=(\beta_1,\dots,\beta_d)^T\in\C^d$, let $H_A(\bbeta)$ be the left ideal of the Weyl algebra $D=\C\langle \bx,\partial_{\bx}\rangle=\C\langle x_1,\dots,x_n,\partial_1,\dots,\partial_n\rangle$ generated by $I_A$ and the Euler operators
\begin{equation}
\label{eq:EulerOperators}
\sum_{j=1}^n a_{ij}\theta_j-\beta_i
\quad
(i=1,\dots,d),
\end{equation}
where $\theta_j=x_j\partial_j$.
The left $D$-module $M_A(\bbeta)=D/H_A(\bbeta)$ is the $A$-hypergeometric system with parameter $\bbeta$.
A formal power series annihilated by $H_A(\bbeta)$ will be called an $A$-hypergeometric series with parameter $\bbeta$.

The systematic study of $A$-hypergeometric systems was initiated by Gel'fand, Graev, Kapranov, and Zelevinsky \cite{GGZ,GZK2,GKZ}.
Logarithm-free series solutions were constructed for generic parameters in \cite{GGZ,GZK2} and in greater generality in \cite{SST}.
Under the homogeneity assumption on $A$, the system $M_A(\bbeta)$ is regular holonomic; see \cite{Hotta,SW}.
The papers \cite{Log,Log2} developed a Frobenius-type perturbation method for constructing logarithmic $A$-hypergeometric series.
We call this construction the intrinsic perturbation construction.
In that construction, one starts from a fake exponent $\bv_0$ and perturbs it within the affine space $\bv_0+L_{\mathbb C}$.
Differentiation with respect to the perturbation parameters then produces polynomials in the logarithms of the variables.

In this paper, we fix a generic weight, a fake exponent $\bv_0$, and an ordered family $N$ of admissible negative supports.
We study the space $C_{\bv_0,N}$ of logarithmic polynomials that can occur as the coefficient of the monomial $\bx^{\bv_0}$ in an $A$-hypergeometric series whose monomial terms satisfy the prescribed support conditions.
Our objective is both to describe this coefficient space explicitly and to determine the subspace realized by the intrinsic perturbation construction, including when this subspace coincides with the full coefficient space.

Our first result reduces the coefficient problem to a finite system.
The Euler operators yield constant-coefficient differential equations for the logarithmic coefficients, while the toric operators relate coefficients attached to different lattice shifts.
Although the resulting system is initially indexed by all elements $\bu\in L$, a suitable normalization makes the normalized coefficient associated with $\bx^{\bv_0+\bu}$ depend only on the negative support of $\bv_0+\bu$.
Since only finitely many negative supports occur, the normalized coefficient equations form a finite system indexed by $N$.

We encode this finite system in a finite free module.
A monomial map from this module to $\C[\partial_{\by}]$ identifies the annihilator of $C_{\bv_0,N}$ with an explicitly defined homogeneous colon ideal.
Equivalently, the coefficient space is the orthogonal complement, or the polynomial solution space, of this colon ideal.
Thus the space of logarithmic polynomials that can occur as the coefficient of $\bx^{\bv_0}$ is determined by an explicitly defined colon ideal constructed from the finite family $N$.
These assertions are stated in \Cref{thm:ann-CvN}.

For the negative support family $N=\cN_{\bv_0}$ determined by the direction $\bw$, we relate this description to the indicial ideal of $H_A(\bbeta)$.
We extract the primary component of the indicial ideal supported at $\bv_0$, shift it to the origin, and replace the Euler variables by constant-coefficient differential operators in the logarithmic variables.
We prove that the orthogonal complement of the resulting Artinian ideal is precisely $C_{\bv_0,\cN_{\bv_0}}$.
Taking annihilators then identifies the shifted indicial component with the colon ideal appearing in \Cref{thm:ann-CvN}.
Thus the finite normalized coefficient system gives an explicit colon-ideal presentation of the local indicial equations at the fake exponent $\bv_0$.
These assertions are stated in \Cref{thm:indicial-ideal-coefficient-space}.

Our second result gives a perturbation-theoretic realization of the full coefficient space.
The intrinsic perturbation construction varies the fake exponent only within $\bv_0+L_{\mathbb C}$.
We instead introduce an ambient perturbation in which the fake exponent is allowed to vary in the full affine space $\bv_0+\mathbb C^n$.
In this formulation, the negative support conditions are encoded by monomial ideals, whereas the Euler equations are encoded by a linear ideal.
Using formal series star-duality, we prove that the resulting construction produces $A$-hypergeometric series and realizes every element of $C_{\bv_0,N}$.
These assertions are stated in \Cref{thm:ambient-perturbation-general-N}.

Finally, we compare the ambient and intrinsic perturbation constructions.
The intrinsic perturbation construction produces a subspace of the coefficient space realized by the ambient perturbation construction and hence a subspace of $C_{\bv_0,N}$.
We prove that the intrinsic perturbation construction realizes the full coefficient space if and only if the ambient colon ideal restricts to the intrinsic colon ideal:
\[
\Phi(\widehat Q_N:\widehat m_N(\bt))=P_N:m_N(\bs).
\]
This assertion is stated in \Cref{prop:L-perturb-general-N}.
Hence the completeness of the intrinsic perturbation construction is reduced to a concrete equality of homogeneous ideals.

We give several sufficient conditions for this equality in \Cref{sec:comparison-intrinsic-ambient-perturbations}.
The examples in \Cref{sec:examples} show that these conditions are not necessary.
We also leave open the question whether the comparison equality can ever fail, or equivalently, whether there exists an ordered negative support family for which the intrinsic perturbation construction does not realize the full coefficient space.

The paper is organized as follows.
In \Cref{sec:duality}, we establish the ordinary and formal series duality conventions used throughout the paper, including the star-operation and its relation to colon ideals.
In \Cref{sec:negative-support-families}, we introduce the negative support data and recall the intrinsic perturbation construction.
In \Cref{sec:coefficient-systems-general-N}, we derive the finite normalized coefficient system, prove the colon ideal description of $C_{\bv_0,N}$, and identify the resulting colon ideal for $N=\cN_{\bv_0}$ with the component of the indicial ideal supported at $\bv_0$ and shifted to the origin.
In \Cref{sec:ambient-perturbation-general-N}, we construct the ambient perturbation series and prove that it realizes $C_{\bv_0,N}$.
In \Cref{sec:comparison-intrinsic-ambient-perturbations}, we compare the intrinsic and ambient perturbation constructions and give sufficient conditions for their coefficient spaces to coincide.
Finally, \Cref{sec:examples} contains two examples, and the negative support computations used there are recorded in the appendix.

\section{Duality}
\label{sec:duality}

In this section, we summarize the duality arising from the pairing between constant-coefficient differential operators and polynomials or formal power series.

\subsection{Pairings}
\label{subsec:duality-conventions-z}

Let $\bz:=(z_1,\dots,z_r)$ be general variables.
We write $\C[\bz]:=\C[z_1,\dots,z_r]$, $\C[[\bz]]:=\C[[z_1,\dots,z_r]]$, and $\C[\partial_{\bz}]:=\C[\partial_{z_1},\dots,\partial_{z_r}]$, where $\partial_{z_j}:=\frac{\partial}{\partial z_j}$ for $j=1,\dots,r$.

We use the usual total degree gradings
\[
\C[\bz]=\bigoplus_{m\ge0}\C[\bz]_m,\quad
\C[\partial_{\bz}]=\bigoplus_{m\ge0}\C[\partial_{\bz}]_m.
\]
The formal power series ring $\C[[\bz]]$ is regarded as the degree-wise direct product
\[
\C[[\bz]]=\prod_{m\ge0}\C[\bz]_m.
\]

We consider the following pairing 
\begin{equation}
\label{eq:z-pairing}
(\,\cdot\,,\,\cdot\,)_{\bz}:\C[\partial_{\bz}]\times\C[[\bz]]\longrightarrow\C,\quad (q,p)_{\bz}:=\left(q(\partial_{\bz})\bullet p(\bz)\right)_{|\bz=\bzero}.
\end{equation}
Here $\bullet$ denotes the usual differential action.
For $\aalpha,\bbeta\in\N^r$, we have 
\[
(\partial_{\bz}^{\aalpha},\bz^{\bbeta})_{\bz}=\delta_{\aalpha,\bbeta}\,\aalpha!.
\]
Here $\partial_{\bz}^{\aalpha}$ and $\bz^{\bbeta}$ denote the usual multi-index notation, $\delta_{\aalpha,\bbeta}$ denotes the Kronecker delta, and $\aalpha!:=\prod_{j=1}^r\alpha_j !$.

Hence the pairing \eqref{eq:z-pairing} is non-degenerate, as is immediate from the monomial bases.
More precisely, for every $m\geq0$, its restriction
\[
\C[\partial_{\bz}]_m\times\C[\bz]_m\longrightarrow\C
\]
is a perfect pairing, whereas
\begin{equation}
\label{eq:z-degree-orthogonality}
(\C[\partial_{\bz}]_m,\C[\bz]_n)_{\bz}=0\quad (m\neq n).
\end{equation}

We shall use the following notation for orthogonal complements.
For $P\subset\C[[\bz]]$ and $Q\subset\C[\partial_{\bz}]$, define
\[
P^{\perp,\C[\partial_{\bz}]}:=\{q\in\C[\partial_{\bz}]\mid(q,p)_{\bz}=0\text{ for every }p\in P\},
\]
\[
Q^{\perp,\C[[\bz]]}:=\{p\in\C[[\bz]]\mid(q,p)_{\bz}=0\text{ for every }q\in Q\}.
\]
If $P\subset\C[\bz]$ or if the orthogonal complement is taken in $\C[\bz]$, we use the analogous notation $P^{\perp,\C[\partial_{\bz}]}$ and $Q^{\perp,\C[\bz]}$.

\subsection{Finite free module duality}
\label{subsec:finite-free-duality}
\label{subsec:polynomial-duality-z}

Let $\Lambda$ be a finite set, and put
\[
 F:=\bigoplus_{\lambda\in\Lambda}\C[\partial_{\bz}]e_\lambda^\ast,
 \quad
 C:=\bigoplus_{\lambda\in\Lambda}\C[\bz]e_\lambda.
\]
The module $C$ is endowed with the usual componentwise differential action.
Thus, for $f=\sum_\lambda f_\lambda e_\lambda^\ast\in F$ and $c=\sum_\lambda c_\lambda e_\lambda\in C$, we set
\[
f\bullet c:=\sum_{\lambda\in\Lambda}f_\lambda\bullet c_\lambda\in\C[\bz].
\]
We also use the componentwise pairing
\[
(f,c)_{\bz}:=\sum_{\lambda\in\Lambda}(f_\lambda,c_\lambda)_{\bz}.
\]
For a submodule $U\subset F$, let
\[
\Sol_C^\bullet(U):=\{c\in C\mid u\bullet c=0\text{ for every }u\in U\}.
\]
For linear subspaces $U\subset F$ and $V\subset C$, let
\[
U^{\perp,C}:=\{c\in C\mid(u,c)_{\bz}=0\text{ for all }u\in U\},
\]
\[
V^{\perp,F}:=\{f\in F\mid(f,v)_{\bz}=0\text{ for all }v\in V\}.
\]
Assign integers $d_\lambda:=\deg(e_\lambda)=\deg(e_\lambda^\ast)$, and give $F$ and $C$ the corresponding shifted gradings
\[
F_m=\bigoplus_{\lambda\in\Lambda}\C[\partial_{\bz}]_{m-d_\lambda}e_\lambda^\ast, \quad C_m=\bigoplus_{\lambda\in\Lambda}\C[\bz]_{m-d_\lambda}e_\lambda.
\]
By the degreewise perfectness of the scalar pairing, the componentwise pairing $F\times C\to\C$ is non-degenerate; indeed, $F_m\times C_m\to\C$ is perfect for every $m\in\Z$, and
\[
(F_m,C_{m'})_{\bz}=0\quad(m\neq m').
\]

For $\lambda_0\in\Lambda$, let $\operatorname{pr}_{\lambda_0}:C\to\C[\bz]$ be the projection onto the $\lambda_0$-component.
For a linear subspace $W\subset\C[\bz]$, put
\[
\Ann_{\C[\partial_{\bz}]}^\bullet(W):=\{f\in\C[\partial_{\bz}]\mid f\bullet w=0 \text{ for all }w\in W\}.
\]

\begin{lemma}
\label{lem:finite-free-duality}
Let $U\subset F$ be a homogeneous $\C[\partial_{\bz}]$-submodule.
Then the following hold.
\begin{enumerate}[(1)]
\item $U^{\perp,C}=\Sol_C^\bullet(U)$.
\item $(U^{\perp,C})^{\perp,F}=U$.
\item For any $\lambda_0\in\Lambda$, 
\begin{equation}
\label{eq:ordinary-component-annihilator}
\Ann_{\C[\partial_{\bz}]}^\bullet\bigl(\operatorname{pr}_{\lambda_0}(U^{\perp,C})\bigr)=\{f\in\C[\partial_{\bz}]\mid fe_{\lambda_0}^\ast\in U\}.
\end{equation}
\end{enumerate}
\end{lemma}

\begin{proof}
Let $u\in U$, $c\in C$, and $q\in\C[\partial_{\bz}]$.
The usual differential action gives
\begin{equation}
\label{eq:ordinary-adjointness}  
(q,u\bullet c)_{\bz}=(qu,c)_{\bz}.
\end{equation}

(1) Suppose that $c\in U^{\perp,C}$.
As $qu\in U$, by \eqref{eq:ordinary-adjointness}, we have $(q,u\bullet c)_\bz=0$ for every $q$.
Non-degeneracy of the pairing on $\C[\partial_{\bz}]\times\C[\bz]$ gives $u\bullet c=0$.
If $u\bullet c=0$, then $(u,c)_\bz=0$ by taking $q=1$ in \eqref{eq:ordinary-adjointness}.
Hence we have the assertion.

(2) Write $U=\bigoplus_m U_m$, where $U_m=U\cap F_m$.
Degree orthogonality and the perfect pairing $F_m\times C_m\to\C$ give
\[
U^{\perp,C}=\bigoplus_m U_m^\perp,\quad (U_m^\perp)^\perp=U_m.
\]
Taking the direct sum over $m$ proves the equality.

(3) If $fe_{\lambda_0}^\ast\in U$, then (1) implies $f\bullet\operatorname{pr}_{\lambda_0}(c)=0$ for all $c\in U^{\perp,C}$.
Conversely, if $f\in \Ann_{\C[\partial_{\bz}]}^\bullet\bigl(\operatorname{pr}_{\lambda_0}(U^{\perp,C})\bigr)$, then 
\[
(fe_{\lambda_0}^\ast,c)_{\bz}=(f,\mathrm{pr}_{\lambda_0}(c))_{\bz}=0
\] 
for all $c\in U^{\perp,C}$.
By (2), we have $fe_{\lambda_0}^\ast\in (U^{\perp,C})^{\perp,F}=U$.
\end{proof}

We now record the following one-component consequences.
\begin{corollary}
\label{cor:one-component-ideal-duality}
Let $I\subset\C[\partial_{\bz}]$ be a homogeneous ideal.
Then
\begin{equation}
\label{eq:ideal-perp-sol}
\Sol_{\C[\bz]}^\bullet(I)=I^{\perp,\C[\bz]},
\end{equation}
and
\begin{equation}
\label{eq:double-perp-ideal}
(I^{\perp,\C[\bz]})^{\perp,\C[\partial_{\bz}]}=I.
\end{equation}
\end{corollary}

\begin{proof}
Apply \Cref{lem:finite-free-duality} to $F=\C[\partial_{\bz}]e^\ast$, $C=\C[\bz]e$, and $U=Ie^\ast$.
\end{proof}

\begin{lemma}
\label{lem:artinian-ideal-orthogonal-duality}
Let $I\subset\C[\partial_{\bz}]$ be an Artinian ideal.
Then $I^{\perp,\C[\bz]}$ is finite-dimensional, and the pairing $(\,\cdot\,,\,\cdot\,)_{\bz}$ induces a perfect pairing
\begin{equation}
\label{eq:artinian-quotient-orthogonal-pairing}
(\C[\partial_{\bz}]/I) \times I^{\perp,\C[\bz]}\longrightarrow\C.
\end{equation}
Moreover,
\begin{equation}
\label{eq:artinian-ideal-orthogonal-recovery}
\Ann_{\C[\partial_{\bz}]}^\bullet\left(I^{\perp,\C[\bz]}\right)=I.
\end{equation}
\end{lemma}

\begin{proof}
Since $I$ is Artinian, there exists an integer $N\geq1$ such that $\langle\partial_{z_1},\dots,\partial_{z_r}\rangle^N\subset I$.
Put
\[
R_{<N}:=\bigoplus_{m=0}^{N-1}\C[\partial_{\bz}]_m,\quad S_{<N}:=\bigoplus_{m=0}^{N-1}\C[\bz]_m.
\]
The pairing $R_{<N}\times S_{<N}\to\C$ is perfect by the degreewise perfectness of the pairing $(\,\cdot\,,\,\cdot\,)_{\bz}$.

Every residue class in $\C[\partial_{\bz}]/I$ has a representative in $R_{<N}$.
Moreover, if $p\in I^{\perp,\C[\bz]}$, then $\left(\partial_{\bz}^{\aalpha},p\right)_{\bz}=0$ for every $\aalpha\in\N^r$ with $|\aalpha|\geq N$, because $\partial_{\bz}^{\aalpha}\in I$.
Hence $p\in S_{<N}$, and therefore $I^{\perp,\C[\bz]}$ is finite-dimensional.

Let $I_{<N}:=I\cap R_{<N}$.
Since every term of degree at least $N$ belongs to $\langle\partial_{z_1},\dots,\partial_{z_r}\rangle^N\subset I$, truncation to degrees less than $N$ induces an isomorphism $\C[\partial_{\bz}]/I\simeq R_{<N}/I_{<N}$.
Under the perfect pairing $R_{<N}\times S_{<N}\to\C$, the orthogonal complement of $I_{<N}$ is precisely $I^{\perp,\C[\bz]}$.
Finite-dimensional linear duality therefore shows that \eqref{eq:artinian-quotient-orthogonal-pairing} is a perfect pairing.

We now prove \eqref{eq:artinian-ideal-orthogonal-recovery}.
Let $f\in I$ and $p\in I^{\perp,\C[\bz]}$.
By \Cref{lem:finite-free-duality} (1), $p\in\Sol^{\bullet}_{\C[\bz]}(I)$. 
Hence $f\bullet p=0$.
Thus
\[
I\subset\Ann_{\C[\partial_{\bz}]}^\bullet\left(I^{\perp,\C[\bz]}\right).
\]

Conversely, let $f\notin I$.
Then the residue class of $f$ in $\C[\partial_{\bz}]/I$ is nonzero.
By the perfectness of \eqref{eq:artinian-quotient-orthogonal-pairing}, there exists $p\in I^{\perp,\C[\bz]}$ such that $(f,p)_{\bz}\neq0$.
In particular, $f\bullet p\neq0$, and hence
\[
f\notin\Ann_{\C[\partial_{\bz}]}^\bullet\left(I^{\perp,\C[\bz]}\right).
\]
This proves the reverse inclusion.
\end{proof}

\begin{lemma}
\label{lem:ordinary-ann-equals-perp-for-stable-subspace}
Let $W\subset\C[\bz]$ be a $\C[\partial_{\bz}]$-stable linear subspace.
Then
\begin{equation}
\label{eq:ordinary-ann-equals-perp-for-stable-subspace}
\Ann_{\C[\partial_{\bz}]}^\bullet(W)=W^{\perp,\C[\partial_{\bz}]}.
\end{equation}
\end{lemma}

\begin{proof}
The inclusion from left to right is immediate.
Conversely, let $f\in W^{\perp,\C[\partial_{\bz}]}$ and $w\in W$.
Since $W$ is stable, for every $q\in\C[\partial_{\bz}]$, we have
\[
(q,f\bullet w)_{\bz}=(qf,w)_{\bz}=(f,q\bullet w)_{\bz}=0.
\]
Non-degeneracy gives $f\bullet w=0$, proving the reverse inclusion.
\end{proof}

\begin{proposition}
\label{prop:double-perp-subspace}
Let $W\subset\C[\bz]$ be a graded $\C$-linear subspace.
Then
\begin{equation}
\label{eq:double-perp-subspace}
(W^{\perp,\C[\partial_{\bz}]})^{\perp,\C[\bz]}=W.
\end{equation}
\end{proposition}

\begin{proof}
This follows from the perfect pairings $\C[\partial_{\bz}]_m\times\C[\bz]_m\to\C$ for all $m$.
\end{proof}

\subsection{Formal series star-duality and colon ideals}
\label{subsec:formal-star-duality-colon}

We recall the star-operation introduced in \cite[Lemmas~3.8 and~3.9]{Log2}.

Let $\xxi=(\xi_1,\dots,\xi_r)$ be auxiliary variables.
For $m\in\C[\bz]$ and $q\in\C[\partial_{\bz}]$, define the star-operation by
\begin{equation}
\label{eq:star-z-action}
m\star q:=\left(m(\partial_{\xxi})\bullet q(\xxi)\right)_{|\xxi=\partial_{\bz}}\in\C[\partial_{\bz}].
\end{equation}
Since $q$ has finite degree, this operation extends to
\[
\C[[\bz]]\times\C[\partial_{\bz}]\longrightarrow\C[\partial_{\bz}].
\]
The extended star-operation is associative in the sense that
\begin{equation}
\label{eq:star-associativity}
(m_1m_2)\star q=m_1\star(m_2\star q)\quad (m_1,m_2\in\C[[\bz]],\ q\in\C[\partial_{\bz}]).
\end{equation}
It also satisfies the adjointness relation
\begin{equation}
\label{eq:starAdjoint}
(m\star q,p)_{\bz}=(q,mp)_{\bz}\quad(m,p\in\C[[\bz]],\ q\in\C[\partial_{\bz}]).
\end{equation}
Thus multiplication in $\C[[\bz]]$ and the star-operation on $\C[\partial_{\bz}]$ are transposes of one another.

For an ideal $P\subset\C[[\bz]]$, define
\[
\Sol_{\C[\partial_{\bz}]}^\star(P):=\{q\in\C[\partial_{\bz}]\mid p\star q=0\text{ for all }p\in P\}.
\]

\begin{lemma}
\label{lem:Pperp}
Let $P\subset\C[[\bz]]$ be an ideal.
Then $P^{\perp,\C[\partial_{\bz}]}$ is a star-submodule of $\C[\partial_{\bz}]$, and
\begin{equation}
\label{eq:Pperp-star-characterization}
P^{\perp,\C[\partial_{\bz}]}=\Sol_{\C[\partial_{\bz}]}^\star(P).
\end{equation}
If $P$ is homogeneous, then $P^{\perp,\C[\partial_{\bz}]}$ is graded.
\end{lemma}

\begin{proof}
For $m\in\C[[\bz]]$, $q\in P^{\perp,\C[\partial_{\bz}]}$, and
$p\in P$, \eqref{eq:starAdjoint} gives
\[
 (m\star q,p)_{\bz}=(q,mp)_{\bz}=0.
\]
Hence $P^{\perp,\C[\partial_{\bz}]}$ is a star-submodule.

If $p\star q=0$ for every $p\in P$, then $(q,p)_{\bz}=(p\star q,1)_{\bz}=0$, so the right-hand side of \eqref{eq:Pperp-star-characterization} is contained in the left-hand side.
Conversely, let $q\in P^{\perp,\C[\partial_{\bz}]}$, $p\in P$, and put $r=p\star q$.
If $r\neq0$, let $d=\deg(r)$.
For every $\aalpha\in\N^r$ with $|\aalpha|=d$, the associativity \eqref{eq:star-associativity} gives
\[
\bz^{\aalpha}\star r=\bz^{\aalpha}\star(p\star q)=(\bz^{\aalpha}p)\star q.
\]
Since $\bz^{\aalpha}p\in P$, the adjointness relation \eqref{eq:starAdjoint} gives
\[
\bz^{\aalpha}\star r=(\bz^{\aalpha}p)\star q=(q,\bz^{\aalpha}p)_{\bz}\cdot1=0.
\]
This contradicts the definition of $d$.
Thus $p\star q=0$.
Finally, if $P$ is homogeneous, gradedness follows immediately from
\eqref{eq:z-degree-orthogonality}.
\end{proof}

\begin{proposition}
\label{prop:Duality}
Let $P\subset\C[[\bz]]$ be a homogeneous ideal, and let
$Q\subset\C[\partial_{\bz}]$ be a graded star-submodule.  Then
\begin{equation}
\label{eq:double-perp-P}
 (P^{\perp,\C[\partial_{\bz}]})^{\perp,\C[[\bz]]}=P,
\end{equation}
and
\begin{equation}
\label{eq:double-perp-Q}
 (Q^{\perp,\C[[\bz]]})^{\perp,\C[\partial_{\bz}]}=Q.
\end{equation}
\end{proposition}

\begin{proof}
Both assertions follow from \eqref{eq:z-degree-orthogonality} and the perfect pairings $\C[\partial_{\bz}]_m\times\C[\bz]_m\to\C$ for all $m$.
\end{proof}

\begin{proposition}
\label{prop:ColonIdeal}
Let $P$ be a homogeneous ideal of $\C[[\bz]]$, and let $m\in\C[[\bz]]$ be homogeneous.
Then
\begin{equation}
\label{eq:colon-star-first}
(m\star P^{\perp,\C[\partial_{\bz}]})^{\perp,\C[[\bz]]}=P:m,
\end{equation}
and
\begin{equation}
\label{eq:colon-star-second}
m\star P^{\perp,\C[\partial_{\bz}]}=(P:m)^{\perp,\C[\partial_{\bz}]}.
\end{equation}
\end{proposition}

\begin{proof}
For $p\in\C[[\bz]]$, the adjointness relation
\eqref{eq:starAdjoint} and \Cref{prop:Duality} give
\begin{align*}
p\in(m\star P^{\perp,\C[\partial_{\bz}]})^{\perp,\C[[\bz]]}&\iff(q,mp)_{\bz}=0 \text{ for every }q\in P^{\perp,\C[\partial_{\bz}]}\\
&\iff mp\in P\\
&\iff p\in P:m.
\end{align*}
This proves \eqref{eq:colon-star-first}.
The star-submodule $m\star P^{\perp,\C[\partial_{\bz}]}$ is graded, and $P:m$ is a homogeneous ideal.
Taking orthogonal complements on both sides of \eqref{eq:colon-star-first} and using \Cref{prop:Duality}, we obtain \eqref{eq:colon-star-second}.
\end{proof}

\subsection{Lattice coordinates and comparison maps}
\label{subsec:t-variables-At-perp-comparison}

Recall that
\[
L=\Ker_{\Z}A\subset\Z^n,\quad L_{\C}=\C\otimes_{\Z}L=\Ker_{\C}A\subset\C^n.
\]
Let $\bt=(t_1,\dots,t_n)^T$, and let $\langle A\bt\rangle$ be the ideal of $\C[[\bt]]$ generated by the linear forms
\begin{equation}
\label{eq:At}  
\ba^{(i)}\bt:=\sum_{j=1}^n a_{ij}t_j\quad(i=1,\dots,d),
\end{equation}
where $\ba^{(i)}$ is the $i$-th row of $A$.
We denote by $\C[L]\subset\C[\partial_{\bt}]$ the symmetric algebra generated by
\[
L_{\C}\cdot\partial_{\bt}:=\left\langle\bu\cdot\partial_{\bt}\,\middle|\,\bu\in L_{\C}\right\rangle_{\C}.
\]
Here $\bu\cdot\partial_{\bt}:=\sum_{j=1}^{n}u_j\partial_{t_j}$ for $\bu=(u_1,\dots,u_n)^T\in\C^n$.
\begin{proposition}
\label{prop:Ay-perp}
One has
\begin{equation}
\label{eq:At-perp}
 \langle A\bt\rangle^{\perp,\C[\partial_{\bt}]}=\C[L].
\end{equation}
\end{proposition}

\begin{proof}
By \Cref{lem:Pperp}, a polynomial $q\in\C[\partial_{\bt}]$ belongs to $\langle A\bt\rangle^{\perp,\C[\partial_{\bt}]}$ if and only if
\[
(\ba^{(i)}\bt)\star q=\left((\ba^{(i)})^T\cdot\partial_{\bt}\right)\bullet q(\bt)=0 \quad(i=1,\dots,d).
\]
By \cite[Lemma~5.1]{LogFree}, this is equivalent to $q\in\C[L]$.
\end{proof}

Fix a basis
\[
 B=\{\bb^{(1)},\dots,\bb^{(n-d)}\}
\]
of $L_{\C}$, and regard $B$ also as the matrix whose columns are these basis vectors.
For $\bs=(s_1,\dots,s_{n-d})^T$, put
\begin{equation}
\label{eq:Bs}
(B\bs)_j:=\sum_{k=1}^{n-d}b_j^{(k)}s_k\quad(j=1,\dots,n).
\end{equation}
Define two $\C$-algebra homomorphisms
\[
\Phi=\Phi_B:\C[[\bt]]\longrightarrow\C[[\bs]],\quad\Phi(f):=f((B\bs)_1,\dots,(B\bs)_n),
\]
and
\[
\Psi=\Psi_B:\C[\partial_{\bs}]\longrightarrow\C[\partial_{\bt}],\quad\Psi(q):=q(\bb^{(1)}\cdot\partial_{\bt},\dots,\bb^{(n-d)}\cdot\partial_{\bt}).
\]
Then $\Phi$ is surjective with $\Ker(\Phi)=\langle A\bt\rangle$, while $\Psi$ is injective with $\operatorname{Im}(\Psi)=\C[L]$.

\begin{proposition}
\label{prop:Ay-perp2}
The following hold.
\begin{enumerate}[(1)]
\item For $f\in\C[[\bt]]$ and $q\in\C[\partial_{\bs}]$,
\begin{equation}
\label{eq:Phi-Psi-adjoint}
 (q,\Phi(f))_{\bs}=(\Psi(q),f)_{\bt}.
\end{equation}
\item If $P\subset\C[[\bt]]$ is an ideal containing $\langle A\bt\rangle$, then
\[
P^{\perp,\C[\partial_{\bt}]}\subset\C[L],
\]
and
\begin{equation}
\label{eq:Phi-Psi-perp}
(\Phi(P))^{\perp,\C[\partial_{\bs}]}=\Psi^{-1}\left(P^{\perp,\C[\partial_{\bt}]}\right).
\end{equation}
Equivalently,
\begin{equation}
\label{eq:Phi-Psi-perp-2}
\Psi\left((\Phi(P))^{\perp,\C[\partial_{\bs}]}\right)=P^{\perp,\C[\partial_{\bt}]}.
\end{equation}
\end{enumerate}
\end{proposition}

\begin{proof}
(1) For $q=\partial_{s_i}$, the chain rule gives
\[
\partial_{s_i}\bullet\Phi(f)=\Phi\left(\left(\sum_{j=1}^n b_j^{(i)}\partial_{t_j}\right)\bullet f\right)=\Phi\bigl(\Psi(\partial_{s_i})\bullet f\bigr).
\]
By linearity and multiplicativity, the same identity holds for every
$q\in\C[\partial_{\bs}]$.
Evaluating at the origin proves the assertion.

(2) If $P\supset\langle A\bt\rangle$, then \Cref{prop:Ay-perp} gives
\[
P^{\perp,\C[\partial_{\bt}]}\subset\langle A\bt\rangle^{\perp,\C[\partial_{\bt}]}=\C[L].
\]
For $q\in\C[\partial_{\bs}]$, \eqref{eq:Phi-Psi-adjoint} gives
\begin{align*}
q\in(\Phi(P))^{\perp,\C[\partial_{\bs}]}&\iff (\Psi(q),p)_{\bt}=0\text{ for all }p\in P\\
&\iff \Psi(q)\in P^{\perp,\C[\partial_{\bt}]}.
\end{align*}
This proves \eqref{eq:Phi-Psi-perp}.
Since $P^{\perp,\C[\partial_{\bt}]}\subset\C[L]=\operatorname{Im}(\Psi)$, we have \eqref{eq:Phi-Psi-perp-2}.
\end{proof}

\section{Negative support families and intrinsic perturbation}
\label{sec:negative-support-families}

In this section, we recall the part of the perturbation construction of \cite{Log,Log2} which will be used in this paper.
The construction is intrinsic in the sense that the fake exponent is perturbed inside the affine space $\bv_0+L_{\C}$ after choosing coordinates on $L_{\C}$.
Thus the variables $\bs$ used below are the coordinates on $L_{\C}$ determined by a fixed basis of $L_{\C}$, as in the comparison map $\Phi_B$ introduced in \Cref{subsec:t-variables-At-perp-comparison}.

The purpose of this section is only to recall the intrinsic $\bs$-side construction in the form needed later.
The corresponding coefficient system in the $\partial_{\by}$-variables will be introduced in \Cref{sec:coefficient-systems-general-N}, and the ambient $\bt$-side construction will be introduced in \Cref{sec:ambient-perturbation-general-N}.
This separation keeps the roles of the three sets of variables distinct: $\bs$ is used for intrinsic perturbation, $\by$ for logarithmic coefficients, and $\bt$ for ambient perturbation.

We begin by recalling the notation concerning fake exponents and negative supports.
The detailed definitions and background are given in \cite{Log,SST}, and the perturbation construction used below is the one developed in \cite{Log,Log2}.

\subsection{Fake exponents and negative supports}
\label{subsec:fake-exponents-negative-supports}

We recall only the notation needed below.
Fix a generic weight $\bw\in\R^n$.
Let $\bttheta=(\theta_1,\dots,\theta_n),\quad \theta_j=x_j\partial_j$.
A vector $\bv\in\C^n$ is called a fake exponent of $H_A(\bbeta)$ with respect to $\bw$ if it is a zero of the fake indicial ideal
\[
\find_{\bw}(H_A(\bbeta))=\bigl(D\,\ini_{\bw}(I_A)\cap\C[\bttheta]\bigr)+\langle A\bttheta-\bbeta\rangle\subset\C[\bttheta].
\]
For details, see \cite{Log,SST}.

For a vector $\bv=(v_1,\dots,v_n)^T\in\C^n$, its negative support is
\[
\nsupp(\bv):=\{j\in\{1,\dots,n\}\mid v_j\in\Z_{<0}\}.
\]

Throughout the rest of the paper, we fix a fake exponent $\bv_0$ of $H_A(\bbeta)$ with respect to the fixed generic weight $\bw$.
For $\bu\in L$, put
\begin{equation}
\label{eq:Iu-definition}
I_{\bu}:=\nsupp(\bv_0+\bu).
\end{equation}
In particular, $I_{\bzero}=\nsupp(\bv_0)$.
We also put
\begin{equation}
\label{eq:NS-definition}
\NS=\NS_{\bv_0}:=\{I_{\bu}\mid \bu\in L\}.
\end{equation}

Following the support restrictions appearing in the perturbation construction of \cite{Log,Log2}, we shall use the following distinguished subfamily of $\NS$:
\begin{equation}
\label{eq:canonical-negative-support-family}
\cN_{\bv_0}:=\left\{I\in\NS\,\middle|\,\bw\cdot\bu\ge0 \text{ for every }\bu\in L\text{ with }I_{\bu}=I\right\}.
\end{equation}
By the argument used in \cite[Proposition~4.4]{Log}, one has $I_{\bzero}\in\cN_{\bv_0}$.

\begin{definition}
\label{def:ordered-negative-support-family}
A subset $N\subset\NS$ is called a negative support family for $\bv_0$ if
\begin{equation}
\label{eq:negative-support-family-condition}
I_{\bzero}\in N\subset\cN_{\bv_0}.
\end{equation}
It is said to be ordered if, for $I\in N$ and $J\in\NS$, the inclusion $J\subset I$ implies $J\in N$.
\end{definition}

Unless otherwise stated, $N$ will denote an ordered negative support family for $\bv_0$.  We put
\begin{equation}
\label{eq:Nc-KN-definition}
N^c:=\NS\setminus N,
\quad
K_N:=\bigcap_{I\in N}I.
\end{equation}

\subsection{The intrinsic \texorpdfstring{$\bs$}{s}-side data}
\label{subsec:intrinsic-s-side-data}

We now attach to the ordered negative support family $N$ the $\bs$-side monomial data used in the intrinsic perturbation construction. 
Fix, once and for all, a basis $B=\{\bb^{(1)},\dots,\bb^{(n-d)}\}$ of $L_{\C}$.
As in \Cref{subsec:t-variables-At-perp-comparison}, we also regard $B$ as the $n\times(n-d)$ matrix whose columns are $\bb^{(1)},\dots,\bb^{(n-d)}$.
Let $\bs=(s_1,\dots,s_{n-d})$ be perturbation variables.
Then $(B\bs)_1,\dots,(B\bs)_{n}$ (see \eqref{eq:Bs}) are the coordinate functions of the intrinsic perturbation $\bv_0+B\bs$ inside the affine space $\bv_0+L_{\C}$.

For a subset $S\subset\{1,\dots,n\}$, we use the notation
\[
(B\bs)^S:=\prod_{j\in S}(B\bs)_j,
\]
with the convention that the empty product is $1$.
We define
\begin{equation}
\label{eq:mNs}
m_N(\bs):=(B\bs)^{I_{\bzero}\setminus K_N}=\prod_{j\in I_{\bzero}\setminus K_N}(B\bs)_j\in\C[[\bs]],
\end{equation}
and
\begin{equation}
\label{eq:PNs}
P_N:=\left\langle(B\bs)^{(I\cup J)\setminus K_N}\,\middle|\,I\in N,\ J\in N^c\right\rangle \subset\C[[\bs]].
\end{equation}

Since each $(B\bs)_j$ is homogeneous of degree $1$, the element $m_N(\bs)$ is homogeneous of degree $|I_{\bzero}\setminus K_N|$.
Moreover, since each $(B\bs)^{(I\cup J)\setminus K_N}$ is homogeneous, the ideal $P_N$ is homogeneous.

\subsection{The intrinsic perturbation series}
\label{subsec:intrinsic-perturbation-series}

For $\bu\in\Z^n$, write
\[
\bu=\bu_+-\bu_-,\quad\bu_+,\bu_-\in\N^n,\quad\supp(\bu_+)\cap\supp(\bu_-)=\emptyset.
\]
For tuples $\aalpha=(\alpha_1,\dots,\alpha_n)\in \C^n$ and $\nnu=(\nu_1,\dots,\nu_n)\in\N^n$, we use the falling factorial notation
\[
[\aalpha]_{\nnu}:=\prod_{j=1}^n\prod_{k=0}^{\nu_j-1}(\alpha_j-k),
\]
where the empty product is understood to be $1$.

For $\bu\in L$, define
\begin{equation}
\label{eq:intrinsic-coefficient-au}
a_{\bu}(\bs):=\frac{[\bv_0+B\bs]_{\bu_-}}{[\bv_0+B\bs+\bu]_{\bu_+}}.
\end{equation}
This is the coefficient used in the intrinsic perturbation construction, and we shall call it the \emph{intrinsic perturbation coefficient}.

Let $N$ be the fixed ordered negative support family for $\bv_0$.
The \emph{intrinsic perturbation series associated with $N$} is
\begin{equation}
\label{eq:intrinsic-perturbation-series}
\widetilde F_N(\bx,\bs):=m_N(\bs)\sum_{\bu\in L,\,I_{\bu}\in N}a_{\bu}(\bs)\,\bx^{\bv_0+B\bs+\bu}.
\end{equation}

Regard $\log\bx=(\log x_1,\dots,\log x_n)$ as a row vector.  Then
\[
(\log\bx)B=\left(\sum_{j=1}^n b_j^{(1)}\log x_j,\dots,\sum_{j=1}^n b_j^{(n-d)}\log x_j\right).
\]
For $r(\partial_{\bs})\in\C[\partial_{\bs}]$, we write $r((\log\bx)B)$ for the polynomial obtained by the substitution
\[
\partial_{s_k}\longmapsto\sum_{j=1}^n b_j^{(k)}\log x_j\quad(k=1,\dots,n-d).
\]
By \cite[Proposition~3.10]{Log2}, for
$q(\partial_{\bs})\in\C[\partial_{\bs}]$, we have
\begin{equation}
\label{eq:intrinsic-leading-coefficient-formula}
\left(q(\partial_{\bs})\bullet\bigl(m_N(\bs)\bx^{\bv_0+B\bs}\bigr)\right)_{|\bs=\bzero}=\bx^{\bv_0}\bigl(m_N(\bs)\star q(\partial_{\bs})\bigr)\bigl((\log\bx)B\bigr).
\end{equation}

\begin{theorem}
\label{thm:intrinsic-perturbation-general-N}
Let $N$ be an ordered negative support family for $\bv_0$.
For every
\[
q(\partial_{\bs})\in P_N^{\perp,\C[\partial_{\bs}]}=\Sol_{\C[\partial_{\bs}]}^\star(P_N),
\]
the series
\[
q(\partial_{\bs})\bullet\widetilde F_N(\bx,\bs)_{|\bs=\bzero}
\]
is an $A$-hypergeometric series in direction $\bw$.
Moreover, the space of coefficients of $\bx^{\bv_0}$ obtained in this way is precisely
\begin{equation}
\label{eq:intrinsic-coefficient-space}
\left\{r((\log\bx)B)\,\middle|\,r(\partial_{\bs})\in (P_N:m_N(\bs))^{\perp,\C[\partial_{\bs}]}\right\}.
\end{equation}
\end{theorem}

\begin{proof}
The first assertion is the perturbation construction of \cite[Remark~6.3]{Log} and \cite[Theorem~2.7]{Log2}, rewritten for the ordered negative support family $N$.
By \Cref{lem:Pperp}, $P_N^{\perp,\C[\partial_{\bs}]}=\Sol_{\C[\partial_{\bs}]}^\star(P_N)$, so every $q\in P_N^{\perp,\C[\partial_{\bs}]}$ gives an $A$-hypergeometric series in direction $\bw$.

It remains to identify the coefficient of $\bx^{\bv_0}$.
Among the summands in the defining series \eqref{eq:intrinsic-perturbation-series}, only the summand indexed by $\bu=\bzero$ can contribute to the coefficient of $\bx^{\bv_0}$ after applying $q(\partial_{\bs})$ and evaluating at $\bs=\bzero$.
Moreover, $a_{\bzero}(\bs)=1$.
By \eqref{eq:intrinsic-leading-coefficient-formula}, the space of coefficients of $\bx^{\bv_0}$ obtained from $q\in P_N^{\perp,\C[\partial_{\bs}]}$ is therefore precisely
\[
\left\{r((\log\bx)B)\,\middle|\,r(\partial_{\bs})\in m_N(\bs)\star P_N^{\perp,\C[\partial_{\bs}]}\right\}.
\]
Since $P_N$ and $m_N(\bs)$ are homogeneous,
\Cref{prop:ColonIdeal} gives
\[
m_N(\bs)\star P_N^{\perp,\C[\partial_{\bs}]}=(P_N:m_N(\bs))^{\perp,\C[\partial_{\bs}]}.
\]
Consequently, the space of coefficients of $\bx^{\bv_0}$ obtained in this way is precisely \eqref{eq:intrinsic-coefficient-space}.
\end{proof}

\section{Coefficient systems for an ordered negative support family}
\label{sec:coefficient-systems-general-N}

In this section, we describe the coefficient equations for a fixed ordered negative support family $N$ for $\bv_0$.
Let $\by=(y_1,\dots,y_n)$ denote auxiliary variables representing $\log\bx$.
Throughout this section, all actions of $\C[\partial_{\by}]$ on $\C[\by]$ are the usual constant-coefficient differential actions, denoted by $\bullet$.
No star-operation is used in the $\by$-variables.

To emphasize that the objects introduced in this section belong to the $\partial_{\by}$-side coefficient system, we use fraktur letters for the corresponding ideals and submodules. 
All orthogonal complements in this section are taken with respect to the $\by$-pairing obtained from \Cref{subsec:duality-conventions-z} by putting $\bz=\by$.

For a subset $S\subset\{1,\dots,n\}$, we use the notation
\[
\partial_{\by}^{S}:=\prod_{j\in S}\partial_{y_j},
\]
with the convention that the empty product is $1$.

We consider a logarithmic series of the form
\begin{equation}
\label{eq:logarithmic-series-N}
\phi_N(\bx)=\sum_{\bu\in L, I_{\bu}\in N}\bx^{\bv_0+\bu}\,r_{\bu}(\log\bx),\quad r_{\bu}(\by)\in\C[\by].
\end{equation}
For convenience, we set $r_{\bu}=0$ if $I_{\bu}\in N^c$.
Let $C_{\bv_0,N}$ denote the space of polynomials that occur as the coefficient $r_{\bzero}(\by)$ of $\bx^{\bv_0}$ in a logarithmic series of the form \eqref{eq:logarithmic-series-N} annihilated by $H_A(\bbeta)$.

\subsection{Coefficient comparison}
\label{subsec:coefficient-comparison}

For $\bu,\bu'\in L$, define
\begin{align}
\label{eq:d-uprime-leftarrow-u}
d_{\bu'\leftarrow\bu}&:=[\partial_{\by}+\bv_0+\bu]_{(\bu-\bu')_+}  \notag\\
&:=\prod_{\nu=1}^n\prod_{\mu=1}^{u_\nu-u'_\nu}\left(\partial_{y_\nu}+(\bv_0)_\nu+u_\nu-\mu+1\right) \in\C[\partial_{\by}],
\end{align}
where the product over $\mu$ is understood to be $1$ if $u_\nu-u'_\nu\le0$.

\begin{lemma}
\label{lem:coefficient-comparison}
Let $\phi_N(\bx)=\sum_{\bu\in L, I_{\bu}\in N}\bx^{\bv_0+\bu}r_{\bu}(\log\bx)$ be as in \eqref{eq:logarithmic-series-N}.
Then $\phi_N(\bx)$ is annihilated by $H_A(\bbeta)$ if and only if the following two conditions hold.
\begin{enumerate}[(1)]
\item For any $\bu\in L$ and any $i=1,\dots,d$,
\begin{equation}
\label{eq:coefficient-euler-equation}
(A\partial_{\by})_i\bullet r_{\bu}=0,
\end{equation}
where $(A\partial_{\by})_i:=\sum_{j=1}^n a_{ij}\partial_{y_j}$.
\item For any $\bu,\bu'\in L$,
\begin{equation}
\label{eq:coefficient-transport-equation}
d_{\bu'\leftarrow\bu}\bullet r_{\bu}-d_{\bu\leftarrow\bu'}\bullet r_{\bu'}=0.
\end{equation}
\end{enumerate}
\end{lemma}

\begin{proof}
We first prove the necessity.
Suppose that $\phi_N(\bx)$ is annihilated by $H_A(\bbeta)$.

For $\nu=1,\dots,n$ and $\bu\in L$, direct computation gives
\[
\theta_{x_\nu}\bullet\bigl(\bx^{\bv_0+\bu}r_{\bu}(\log\bx)\bigr)=\bx^{\bv_0+\bu}\bigl((\partial_{y_\nu}+(\bv_0)_\nu+u_\nu)\bullet r_{\bu}\bigr)(\log\bx).
\]
Since $A(\bv_0+\bu)=\bbeta$, the equation $(A\theta-\bbeta)_i\bullet\phi_N(\bx)=0$ implies \eqref{eq:coefficient-euler-equation}.

Next, let $\bu,\bu'\in L$.
Put $\bu'':=\bu'-\bu$ and write $\bu''=\bu''_+-\bu''_-$.
Consider the toric operator $\partial_{\bx}^{\bu''_+}-\partial_{\bx}^{\bu''_-}$.
The coefficient of the common monomial $\bx^{\bv_0+\bu-\bu''_-}=\bx^{\bv_0+\bu'-\bu''_+}$ in $\left(\partial_{\bx}^{\bu''_+}-\partial_{\bx}^{\bu''_-}\right)\bullet\phi_N(\bx)$ is
\[
(d_{\bu\leftarrow\bu'}\bullet r_{\bu'})(\log\bx)-(d_{\bu'\leftarrow\bu}\bullet r_{\bu})(\log\bx).
\]
Since the toric operator annihilates $\phi_N(\bx)$, this coefficient must be zero.
Hence we have \eqref{eq:coefficient-transport-equation}.

Conversely, \eqref{eq:coefficient-euler-equation} and \eqref{eq:coefficient-transport-equation} force every coefficient of
\[
(A\theta-\bbeta)_i\bullet\phi_N\quad (i=1,\dots,d)\quad \text{and}\quad(\partial_{\bx}^{\ggamma_+}-\partial_{\bx}^{\ggamma_-})\bullet\phi_N\quad (\ggamma\in L)
\]
to vanish.
Hence $H_A(\bbeta)\bullet\phi_N=0$.
\end{proof}

\begin{lemma}
\label{lem:d-factorization}
Let $\bu,\bu',\bu''\in L$. 
Then the following hold.
\begin{enumerate}[(1)]
\item There exists a unique polynomial $\widetilde d_{\bu'\leftarrow\bu}\in\C[\partial_{\by}]$ with nonzero constant term such that
\begin{equation}
\label{eq:d-factorization}
d_{\bu'\leftarrow\bu}=\widetilde d_{\bu'\leftarrow\bu}\,\partial_{\by}^{I_{\bu'}\setminus I_{\bu}}.
\end{equation}
Moreover, the operator $\widetilde d_{\bu'\leftarrow\bu}\bullet:\C[\by]\longrightarrow\C[\by]$ is a $\C$-linear automorphism.
\item One has the identities:
\begin{equation}
\label{eq:d-cocycle}
d_{\bu\leftarrow\bu''}d_{\bu''\leftarrow\bu'}d_{\bu'\leftarrow\bu}=d_{\bu\leftarrow\bu'}d_{\bu'\leftarrow\bu''}d_{\bu''\leftarrow\bu},
\end{equation}
and
\begin{equation}
\label{eq:dtilde-cocycle}
\widetilde d_{\bu\leftarrow\bu''}\widetilde d_{\bu''\leftarrow\bu'}\widetilde d_{\bu'\leftarrow\bu}=\widetilde d_{\bu\leftarrow\bu'}\widetilde d_{\bu'\leftarrow\bu''}\widetilde d_{\bu''\leftarrow\bu}.
\end{equation}
\end{enumerate}
\end{lemma}

\begin{proof}
(1) We argue coordinatewise.
Fix $\nu\in\{1,\dots,n\}$.
The factor $\partial_{y_\nu}$ occurs in $\prod_{\mu=1}^{u_\nu-u'_\nu}\left(\partial_{y_\nu}+(\bv_0+\bu)_\nu-\mu+1\right)$ if and only if there exists $\mu$ with $1\le \mu\le u_\nu-u'_\nu$ such that $(\bv_0+\bu)_\nu-\mu+1=0$.
This is equivalent to $(\bv_0+\bu)_\nu\in\N$ and $(\bv_0+\bu')_\nu\in\Z_{<0}$, that is, $\nu\in I_{\bu'}\setminus I_{\bu}$.
Moreover, such a factor occurs at most once for each coordinate $\nu$.
After removing all these monomial factors, the remaining polynomial has nonzero constant term.  This proves the asserted factorization and its uniqueness.

The operator $\widetilde d_{\bu'\leftarrow\bu}\bullet$ is invertible on $\C[\by]$ because its constant term is nonzero and its nonconstant part strictly lowers total degree.
This proves the automorphism assertion.

(2) Fix $\nu$, and put
\[
a=(\bv_0+\bu)_\nu,\quad b=(\bv_0+\bu')_\nu,\quad c=(\bv_0+\bu'')_\nu.
\]
The three differences $a-b$, $b-c$, and $c-a$ are integers whose sum is zero.
Hence both sides of \eqref{eq:d-cocycle} have the same multiset of affine factors in $\partial_{y_\nu}$.
Taking the product over all coordinates gives \eqref{eq:d-cocycle}.
Since the monomial factors on both sides are the same, cancelling them gives \eqref{eq:dtilde-cocycle}.
\end{proof}

\subsection{Normalized coefficient system}
\label{subsec:normalized-coefficient-system}

By \Cref{lem:d-factorization}, for every $\bu\in L$, we define the
normalized coefficient $c_{\bu}$ by
\begin{equation}
\label{eq:normalized-coefficient-definition}
c_{\bu}:=\left(\widetilde d_{\bu\leftarrow\bzero}\bullet\right)^{-1}\left(\widetilde d_{\bzero\leftarrow\bu}\bullet r_{\bu}\right).
\end{equation}
Equivalently,
\begin{equation}
\label{eq:ru-in-terms-of-cu}
r_{\bu}=\left(\widetilde d_{\bzero\leftarrow\bu}\bullet\right)^{-1}\left(\widetilde d_{\bu\leftarrow\bzero}\bullet c_{\bu}\right).
\end{equation}
Since $r_{\bu}=0$ whenever $I_{\bu}\in N^c$, we have
\begin{equation}
\label{eq:cu-zero-outside-N}
c_{\bu}=0\quad(I_{\bu}\in N^c).
\end{equation}
Moreover, we have $c_{\bzero}=r_{\bzero}$.

\begin{proposition}
\label{prop:normalized-coefficient-system}
Let $c_{\bu}$ be the normalized coefficients defined by \eqref{eq:normalized-coefficient-definition}.
The logarithmic series $\phi_N(\bx)$ in \eqref{eq:logarithmic-series-N} is annihilated by $H_A(\bbeta)$ if and only if the following two conditions hold.

\begin{enumerate}[(1)]
\item For any $\bu,\bu'\in L$, the equality $I_{\bu}=I_{\bu'}$ implies $c_{\bu}=c_{\bu'}$.
Thus the normalized coefficient depends only on the negative support.
We shall write $c_I:=c_{\bu}$ if $I=I_{\bu}\in\NS$.
In particular,
\begin{equation}
\label{eq:cI-zero-outside-N}
c_I=0\quad (I\in N^c).
\end{equation}

\item The family $(c_I)_{I\in\NS}$ satisfies
\begin{align}
\label{eq:normalized-euler-equation}
(A\partial_{\by})_i\bullet c_I&=0 \quad (I\in\NS,\ i=1,\dots,d),
\\
\label{eq:normalized-transport-equation}
\partial_{\by}^{I'\setminus I}\bullet c_I-\partial_{\by}^{I\setminus I'}\bullet c_{I'}&=0\quad (I,I'\in\NS).
\end{align}
\end{enumerate}
\end{proposition}

\begin{proof}
Suppose that $\phi_N(\bx)$ is annihilated by $H_A(\bbeta)$.
By \Cref{lem:coefficient-comparison}, the coefficients $r_{\bu}$ satisfy \eqref{eq:coefficient-euler-equation} and \eqref{eq:coefficient-transport-equation}.

Fix $\bu,\bu'\in L$.
By \Cref{lem:d-factorization}, we have
\[
d_{\bu'\leftarrow\bu}=\widetilde d_{\bu'\leftarrow\bu}\partial_{\by}^{I_{\bu'}\setminus I_{\bu}},\quad d_{\bu\leftarrow\bu'}=\widetilde d_{\bu\leftarrow\bu'}\partial_{\by}^{I_{\bu}\setminus I_{\bu'}}.
\]
Substituting \eqref{eq:ru-in-terms-of-cu} into \eqref{eq:coefficient-transport-equation} and using \eqref{eq:dtilde-cocycle}, we obtain
\begin{equation}
\label{eq:normalized-transport-with-unit}
h_{\bu,\bu'}\bullet\left(\partial_{\by}^{I_{\bu'}\setminus I_{\bu}}\bullet c_{\bu}-\partial_{\by}^{I_{\bu}\setminus I_{\bu'}}\bullet c_{\bu'}\right)=0,
\end{equation}
where
\[
h_{\bu,\bu'}:=\widetilde d_{\bzero\leftarrow\bu'}\widetilde d_{\bu'\leftarrow\bu}\widetilde d_{\bu\leftarrow\bzero}=\widetilde d_{\bzero\leftarrow\bu}\widetilde d_{\bu\leftarrow\bu'}\widetilde d_{\bu'\leftarrow\bzero}.
\]
Each factor in $h_{\bu,\bu'}$ has a nonzero constant term.
Hence $h_{\bu,\bu'}\bullet$ is an automorphism of $\C[\by]$, and \eqref{eq:normalized-transport-with-unit} is equivalent to
\begin{equation}
\label{eq:normalized-transport-u}
\partial_{\by}^{I_{\bu'}\setminus I_{\bu}}\bullet c_{\bu}-\partial_{\by}^{I_{\bu}\setminus I_{\bu'}}\bullet c_{\bu'}=0.
\end{equation}
In particular, if $I_{\bu}=I_{\bu'}$, then $c_{\bu}=c_{\bu'}$. 
Thus $c_{\bu}$ depends only on $I_{\bu}$, and we may write $c_I:=c_{\bu}$ whenever $I=I_{\bu}$.
Equation \eqref{eq:normalized-transport-u} then becomes \eqref{eq:normalized-transport-equation}.
Moreover, \eqref{eq:cI-zero-outside-N} follows from \eqref{eq:cu-zero-outside-N}.

The normalizing operators in \eqref{eq:normalized-coefficient-definition} are constant-coefficient automorphisms of $\C[\by]$ and commute with each operator $(A\partial_{\by})_i$.
Therefore \eqref{eq:coefficient-euler-equation} is equivalent to \eqref{eq:normalized-euler-equation}.

Conversely, suppose that the conditions in the statement hold.
For each $\bu\in L$, define
\[
r_{\bu}:=\left(\widetilde d_{\bzero\leftarrow\bu}\bullet\right)^{-1}\left(\widetilde d_{\bu\leftarrow\bzero}\bullet c_{I_{\bu}}\right),
\]
where $c_I:=0$ for $I\in N^c$.
Then $r_{\bu}=0$ whenever $I_{\bu}\in N^c$.
Reversing the preceding identities gives \eqref{eq:coefficient-euler-equation} and \eqref{eq:coefficient-transport-equation}.
Hence $\phi_N(\bx)$ is annihilated by $H_A(\bbeta)$ by \Cref{lem:coefficient-comparison}.
\end{proof}

\subsection{The finite normalized module}
\label{subsec:finite-normalized-module}

By \Cref{prop:normalized-coefficient-system}, the normalized coefficients $(c_I)_{I\in \NS}$ satisfy \eqref{eq:normalized-euler-equation} and \eqref{eq:normalized-transport-equation}.
In addition, $c_I=0$ for $I\in N^c$. 
Thus, the unknown part of the system is the finite family $(c_I)_{I\in N}$.  Therefore, the preceding equations reduce to the following \emph{finite} system:
\begingroup
\mathtoolsset{showonlyrefs=false}
\begin{align}
\tag{\ref*{eq:normalized-euler-equation}}
(A\partial_{\by})_i\bullet c_I&=0\quad(I\in N,\ i=1,\dots,d),\\
\label{eq:finite-normalized-transport}
\partial_{\by}^{I'\setminus I}\bullet c_I-\partial_{\by}^{I\setminus I'}\bullet c_{I'}&=0\quad(I,I'\in N),\\
\label{eq:finite-normalized-boundary}
\partial_{\by}^{J\setminus I}\bullet c_I&=0\quad(I\in N,\ J\in N^c).
\end{align}
\endgroup
The last equation is simply \eqref{eq:normalized-transport-equation} with $I'\in N^c$ and $c_{I'}=0$.

Let
\begin{equation}
\label{eq:FN-CN-definition}
F_N:=\bigoplus_{I\in N}\C[\partial_{\by}]e_I^\ast,\quad C_N:=\bigoplus_{I\in N}\C[\by]e_I.
\end{equation}
We give $F_N$ and $C_N$ the same shifted grading determined by
\begin{equation}
\label{eq:FN-shifted-degree}
\deg(e_I)=\deg(e_I^\ast)=|I\setminus K_N|\quad (I\in N).
\end{equation}
Thus
\begin{align*}
(F_N)_m&=\bigoplus_{I\in N}\C[\partial_{\by}]_{m-|I\setminus K_N|}e_I^\ast\quad (m\in\Z),\\
(C_N)_m&=\bigoplus_{I\in N}\C[\by]_{m-|I\setminus K_N|}e_I\quad (m\in\Z).
\end{align*}
Here $\C[\partial_{\by}]_k=0$ and $\C[\by]_k=0$ for $k<0$.
With these shifted gradings, the componentwise pairing $F_N\times C_N\longrightarrow\C$ satisfies
\[
((F_N)_m,(C_N)_{m'})_{\by}=0
\quad
(m\ne m'),
\]
and the induced pairing $(F_N)_m\times (C_N)_m\longrightarrow\C$ is perfect for every $m\in\Z$.
Thus the finite free module duality of \Cref{subsec:finite-free-duality} applies to this shifted graded setting.

For $I\in N$, define
\begin{equation}
\label{eq:P-N-I-definition}
\fkP_{N,I}:=\left\langle\partial_{\by}^{J\setminus I}\,\middle|\,J\in N^c\right\rangle\subset\C[\partial_{\by}].
\end{equation}
Define
\begin{equation}
\label{eq:TN-definition}
\fkT_N:=\left\langle\partial_{\by}^{I'\setminus I}e_I^\ast-\partial_{\by}^{I\setminus I'}e_{I'}^\ast\,\middle|\,I,I'\in N\right\rangle\subset F_N,
\end{equation}
and put
\begin{equation}
\label{eq:UN-definition}
\fkU_N:=\sum_{I\in N}\bigl(\langle A\partial_{\by}\rangle+\fkP_{N,I}\bigr)e_I^\ast+\fkT_N\subset F_N,
\end{equation}
where
\begin{equation}
\label{eq:Apartial-y-ideal}
\langle A\partial_{\by}\rangle:=\left\langle(A\partial_{\by})_i\,\middle|\,i=1,\dots,d\right\rangle\subset\C[\partial_{\by}].
\end{equation}

The submodule $\fkU_N$ is homogeneous with respect to the shifted grading on $F_N$. 
Indeed, the Euler generators $(A\partial_{\by})_ie_I^\ast$ are homogeneous, the generators $\partial_{\by}^{J\setminus I}e_I^\ast$ are homogeneous, and the two terms of
\[
\partial_{\by}^{I'\setminus I}e_I^\ast-\partial_{\by}^{I\setminus I'}e_{I'}^\ast
\]
have the same shifted degree, because $K_N\subset I\cap I'$.

By construction, $c=\sum_{I\in N}c_Ie_I\in C_N$ satisfies the finite normalized system \eqref{eq:normalized-euler-equation}, \eqref{eq:finite-normalized-transport}, and \eqref{eq:finite-normalized-boundary}
if and only if
\[
c\in\Sol_{C_N}^{\bullet}(\fkU_N).
\]
Since $\fkU_N$ is homogeneous, \Cref{lem:finite-free-duality}
gives
\begin{equation}
\label{eq:UN-solution-orthogonal}
\Sol_{C_N}^{\bullet}(\fkU_N)=\fkU_N^{\perp,C_N}.
\end{equation}

By \Cref{prop:normalized-coefficient-system}, $c=\sum_{I\in N}c_Ie_I$ belongs to $\fkU_N^{\perp,C_N}$ if and only if its corresponding logarithmic series $\phi_N(\bx)$ is annihilated by $H_A(\bbeta)$.

Recall that $C_{\bv_0,N}$ denotes the space of all polynomials which occur as the coefficient $r_{\bzero}(\by)$ of $\bx^{\bv_0}$ in a series $\phi_N(\bx)$ annihilated by $H_A(\bbeta)$.
Since $I_{\bzero}\in N$ and the normalization satisfies $r_{\bzero}=c_{\bzero}$, we obtain
\begin{equation}
\label{eq:CvN-as-projection-before-prop}
C_{\bv_0,N}=\operatorname{pr}_{I_{\bzero}}\left(\fkU_N^{\perp,C_N}\right),
\end{equation}
where $\operatorname{pr}_{I_{\bzero}}:C_N\to\C[\by]$ denotes the projection onto the $I_{\bzero}$-component.

\begin{proposition}
\label{prop:ann-CvN-module}
We have the equality
\begin{equation}
\label{eq:ann-CvN-module}
\Ann_{\C[\partial_{\by}]}^\bullet(C_{\bv_0,N})=\left\{f\in\C[\partial_{\by}]\,\middle|\,f e_{I_{\bzero}}^\ast\in\fkU_N\right\}.
\end{equation}
Equivalently,
\begin{equation}
\label{eq:ann-CvN-module-intersection}
\Ann_{\C[\partial_{\by}]}^\bullet(C_{\bv_0,N})e_{I_{\bzero}}^\ast=\fkU_N\cap\C[\partial_{\by}]e_{I_{\bzero}}^\ast.
\end{equation}
\end{proposition}

\begin{proof}
This follows from $C_{\bv_0,N}=\operatorname{pr}_{I_{\bzero}}(\fkU_N^{\perp,C_N})$ and \Cref{lem:finite-free-duality}.
\end{proof}

\subsection{The monomial elimination map}
\label{subsec:monomial-elimination-map}

We now eliminate the redundant $N$-components in the finite normalized module.
Define a $\C[\partial_{\by}]$-linear map
\begin{equation}
\label{eq:muN-definition}
\mu_N:F_N\longrightarrow \C[\partial_{\by}],\quad \mu_N(e_I^\ast)=\partial_{\by}^{I\setminus K_N} \quad (I\in N).
\end{equation}
With respect to the shifted grading on $F_N$ defined in \eqref{eq:FN-shifted-degree} and the usual total degree grading on $\C[\partial_{\by}]$, the map $\mu_N$ is homogeneous of degree $0$.
Indeed,
\[
\mu_N\left(\C[\partial_{\by}]_{m-|I\setminus K_N|}e_I^\ast\right)\subset\C[\partial_{\by}]_m.
\]

\begin{lemma}
\label{lem:kernel-muN}
We have
\begin{equation}
\label{eq:kernel-muN}
\Ker(\mu_N)=\fkT_N.
\end{equation}
\end{lemma}

\begin{proof}
Put $m_I:=\partial_{\by}^{I\setminus K_N}$ for $I\in N$.
Then $\mu_N(e_I^\ast)=m_I$, and
\[
\operatorname{Im}(\mu_N)=\langle m_I\mid I\in N\rangle.
\]
Thus $\Ker(\mu_N)$ is the syzygy module of the monomial generators $\{m_I\,|\,I\in N\}$, and hence is generated by the pairwise monomial syzygies
\[
\frac{\operatorname{lcm}(m_I,m_{I'})}{m_I}e_I^\ast-\frac{\operatorname{lcm}(m_I,m_{I'})}{m_{I'}}e_{I'}^\ast
\quad(I,I'\in N).
\]
Indeed, after choosing a monomial order, the largest monomial occurring in a relation $\sum_I f_I m_I=0$ can be cancelled by a linear combination of these pairwise syzygies.
Induction on the largest monomial proves the assertion.

Since $K_N\subset I\cap I'$, we have
\[
\operatorname{lcm}(m_I,m_{I'})=\partial_{\by}^{(I\cup I')\setminus K_N},
\]
and hence
\[
\frac{\operatorname{lcm}(m_I,m_{I'})}{m_I}=\partial_{\by}^{I'\setminus I},\quad\frac{\operatorname{lcm}(m_I,m_{I'})}{m_{I'}}=\partial_{\by}^{I\setminus I'}.
\]
Therefore $\Ker(\mu_N)$ is generated by
\[
\partial_{\by}^{I'\setminus I}e_I^\ast-\partial_{\by}^{I\setminus I'}e_{I'}^\ast\quad(I,I'\in N),
\]
which are precisely the generators of $\fkT_N$.
Thus we have the assertion.
\end{proof}

\subsection{The colon formula}
\label{subsec:coefficient-colon-formula}

We now pass from the finite normalized module to an ideal in the polynomial ring $\C[\partial_{\by}]$.
Define
\begin{equation}
\label{eq:frakmN-y-definition}
\fkm_N:=\partial_{\by}^{I_{\bzero}\setminus K_N}\in\C[\partial_{\by}],
\end{equation}
\begin{equation}
\label{eq:frakPN-y-definition}
\fkP_N:=\left\langle\partial_{\by}^{(I\cup J)\setminus K_N}\,\middle|\,I\in N,\ J\in N^c\right\rangle\subset\C[\partial_{\by}],
\end{equation}
and
\begin{equation}
\label{eq:frakMN-y-definition}
\fkM_N:=\left\langle\partial_{\by}^{I\setminus K_N}\,\middle|\,I\in N\right\rangle\subset\C[\partial_{\by}].
\end{equation}
Finally, put
\begin{equation}
\label{eq:frakQN-definition}
\fkQ_N:=\langle A\partial_{\by}\rangle\fkM_N+\fkP_N\subset\C[\partial_{\by}].
\end{equation}
All these ideals are homogeneous.
Indeed, $\fkm_N$ and the generators of $\fkM_N$ and $\fkP_N$ are monomials, while $\langle A\partial_{\by}\rangle$ is generated by linear forms.

\begin{lemma}
\label{lem:muN-UN-QN}
We have
\begin{equation}
\label{eq:muN-UN-QN}
\mu_N(\fkU_N)=\fkQ_N,
\end{equation}
and
\begin{equation}
\label{eq:muN-preimage-QN}
\mu_N^{-1}(\fkQ_N)=\fkU_N.
\end{equation}
\end{lemma}

\begin{proof}
We prove \eqref{eq:muN-UN-QN}.
First, for $I\in N$ and $i=1,\dots,d$, the elements $\mu_N((A\partial_{\by})_i e_I^\ast)=(A\partial_{\by})_i\,\partial_{\by}^{I\setminus K_N}$ generate $\langle A\partial_{\by}\rangle\fkM_N$.

Next, for $I\in N$ and $J\in N^c$, the elements $\mu_N(\partial_{\by}^{J\setminus I}e_I^\ast)=\partial_{\by}^{(I\cup J)\setminus K_N}$ generate $\fkP_N$.  
Finally, since $\mu_N(\fkT_N)=0$ by \Cref{lem:kernel-muN}, we have
\[
\mu_N(\fkU_N)=\langle A\partial_{\by}\rangle\fkM_N+\fkP_N=\fkQ_N.
\]

We prove \eqref{eq:muN-preimage-QN}.
The inclusion $\fkU_N\subset \mu_N^{-1}(\fkQ_N)$ is clear by \eqref{eq:muN-UN-QN}.  
Conversely, let $\eta\in \mu_N^{-1}(\fkQ_N)$.
By \eqref{eq:muN-UN-QN}, there exists $\eta'\in\fkU_N$ such that
\[
\mu_N(\eta')=\mu_N(\eta).
\]
Since $\eta-\eta'\in\Ker(\mu_N)=\fkT_N\subset\fkU_N$ by \Cref{lem:kernel-muN}, we obtain $\eta\in\fkU_N$.
Thus $\mu_N^{-1}(\fkQ_N)\subset\fkU_N$, and the proof is complete.
\end{proof}

\begin{theorem}
\label{thm:ann-CvN}
The annihilator of the coefficient space $C_{\bv_0,N}$ is
\begin{equation}
\label{eq:ann-CvN-colon}
\Ann_{\C[\partial_{\by}]}^\bullet(C_{\bv_0,N})=\fkQ_N:\fkm_N,
\end{equation}
and
\begin{equation}
\label{eq:CvN-solution-colon}
C_{\bv_0,N}=(\fkQ_N:\fkm_N)^{\perp,\C[\by]}=\Sol_{\C[\by]}^\bullet(\fkQ_N:\fkm_N).
\end{equation}
\end{theorem}

\begin{proof}
Let $f\in\C[\partial_{\by}]$. 
Using \Cref{lem:muN-UN-QN}, we obtain
\begin{align}
f e_{I_{\bzero}}^\ast\in\fkU_N&\iff f\fkm_N=\mu_N(f e_{I_{\bzero}}^\ast)\in\fkQ_N\notag\\
&\iff f\in\fkQ_N:\fkm_N.
\label{eq:ann-CvN-colon-proof}
\end{align}
By \Cref{prop:ann-CvN-module}, we have 
\[
\Ann_{\C[\partial_{\by}]}^\bullet(C_{\bv_0,N})=\left\{f\in\C[\partial_{\by}]\,\middle|\,f e_{I_{\bzero}}^\ast\in\fkU_N\right\}=\fkQ_N:\fkm_N.
\]

It remains to prove \eqref{eq:CvN-solution-colon}.
By \eqref{eq:CvN-as-projection-before-prop} and \eqref{eq:UN-solution-orthogonal}, we have
\[
C_{\bv_0,N}=\operatorname{pr}_{I_{\bzero}}\left(\Sol_{C_N}^{\bullet}(\fkU_N)\right)=\operatorname{pr}_{I_{\bzero}}\left(\fkU_N^{\perp,C_N}\right).
\]
Since $\fkU_N$ is homogeneous, its orthogonal complement $\fkU_N^{\perp,C_N}$ is graded.
Moreover, $\Sol_{C_N}^{\bullet}(\fkU_N)$ is stable under the componentwise $\C[\partial_{\by}]$-action because constant-coefficient differential operators commute.
Since the projection $\operatorname{pr}_{I_{\bzero}}:C_N\to\C[\by]$ is homogeneous and $\C[\partial_{\by}]$-linear, it follows that $C_{\bv_0,N}$ is a graded $\C[\partial_{\by}]$-submodule of $\C[\by]$.

Therefore, by \Cref{lem:ordinary-ann-equals-perp-for-stable-subspace} and \eqref{eq:ann-CvN-colon}, we have
\[
(C_{\bv_0,N})^{\perp,\C[\partial_{\by}]}=\Ann_{\C[\partial_{\by}]}^\bullet(C_{\bv_0,N})=\fkQ_N:\fkm_N.
\]
Taking orthogonal complements and applying \Cref{prop:double-perp-subspace}, we obtain
\[
C_{\bv_0,N}=(\fkQ_N:\fkm_N)^{\perp,\C[\by]}.
\]
Finally, \Cref{cor:one-component-ideal-duality} gives
\[
(\fkQ_N:\fkm_N)^{\perp,\C[\by]}=\Sol_{\C[\by]}^\bullet(\fkQ_N:\fkm_N),
\]
which proves \eqref{eq:CvN-solution-colon}.
\end{proof}

We now relate the colon ideal in \Cref{thm:ann-CvN} to the indicial
ideal of $H_A(\bbeta)$.

Let $\ind_{\bw}(H_A(\bbeta))\subset\C[\bttheta]$ denote the indicial ideal of $H_A(\bbeta)$ with respect to $\bw$.
For a point $\bv\in\C^n$, let
\[
\mathfrak m_{\bv}:=\left\langle\theta_1-v_1,\dots,\theta_n-v_n\right\rangle\subset\C[\bttheta].
\]
We define the component of $\ind_{\bw}(H_A(\bbeta))$ supported at
$\bv_0$ by
\[
\ind_{\bw}(H_A(\bbeta))_{|\bv_0}:=\left(\ind_{\bw}(H_A(\bbeta))\C[\bttheta]_{\mathfrak m_{\bv_0}}\right)\cap\C[\bttheta].
\]

Let
\[
\tau_{\bv_0}:\C[\bttheta]\longrightarrow\C[\bttheta],\quad \tau_{\bv_0}(f)(\bttheta):=f(\bttheta+\bv_0),
\]
and put
\begin{equation}
\label{eq:shifted-local-indicial-ideal}
\ind_{\bw}(H_A(\bbeta))_{|\bv_0\to\bzero}:=\tau_{\bv_0}\left(\ind_{\bw}(H_A(\bbeta))_{|\bv_0}\right).
\end{equation}
Thus $\ind_{\bw}(H_A(\bbeta))_{|\bv_0\to\bzero}$ is supported at the origin.

Finally, if $J$ is an ideal of $\C[\bttheta]$, we write
\[
J_{\bttheta\to\partial_{\by}}\subset\C[\partial_{\by}]
\]
for the ideal obtained from $J$ by the substitution $\theta_j\mapsto\partial_{y_j}$ for $j=1,\dots,n$.

Let 
\begin{equation}
\label{eq:logarithmic-hypergeometric-series-in direction-w}
\phi(\bx)=\sum_{\bu\in L}\bx^{\bv_0+\bu}r_{\bu}(\log\bx)  
\end{equation}
be an $A$-hypergeometric series in direction $\bw$, where $r_{\bu}(\by)=0$ whenever the corresponding term does not occur.

It was shown in \cite[Corollary~4.3]{Log} that, if $r_{\bu}\ne 0$ then every $\bu'\in L$ satisfying \[\nsupp(\bv_0+\bu')=\nsupp(\bv_0+\bu) \] belongs to the semigroup $C(\bw)$.
Here $C(\bw)$ is generated by the lattice differences corresponding to the reduced Gröbner basis of $I_A$ with respect to $\bw$.
Note that the Gr\"{o}bner basis dealt in \cite[Section~4]{Log} should be the reduced one.
The following lemma gives a direct formulation of this property in terms of $\cN_{\bv_0}$ (see \eqref{eq:canonical-negative-support-family}).

\begin{lemma} \label{lem:direction-support-contained-canonical-family}
Let $\phi(\bx)=\sum_{\bu\in L}\bx^{\bv_0+\bu}r_{\bu}(\log\bx)$ be as in \eqref{eq:logarithmic-hypergeometric-series-in direction-w}.
Then $r_{\bu}\neq0$ implies $I_{\bu}\in\cN_{\bv_0}$.
Consequently, every negative support occurring in $\phi(\bx)$ belongs to $\cN_{\bv_0}$. 
\end{lemma}
\begin{proof}
Let $\bu\in L$ satisfy $r_{\bu}\neq 0$.
Suppose, to the contrary, that $I_{\bu}\notin\cN_{\bv_0}$.
By \eqref{eq:canonical-negative-support-family}, there exists $\bu'\in L$ such that $I_{\bu'}=I_{\bu}$ and $\bw\cdot\bu'<0$.
Since $\phi(\bx)$ is a series in direction $\bw$, the term indexed by $\bu'$ does not occur, and hence $r_{\bu'}=0$.

The corresponding normalized coefficient therefore satisfies $c_{\bu'}=0$.
By \Cref{prop:normalized-coefficient-system} (1), the equality $I_{\bu'}=I_{\bu}$ implies $c_{\bu}=c_{\bu'}=0$.
Since the normalizing operators in \eqref{eq:ru-in-terms-of-cu} are automorphisms, it follows that $r_{\bu}=0$, a contradiction.
Therefore $I_{\bu}\in\cN_{\bv_0}$.
\end{proof}




\begin{theorem}
\label{thm:indicial-ideal-coefficient-space}
We have
\begin{equation}
\label{eq:Cv0-indicial-perp}
C_{\bv_0,\cN_{\bv_0}}=(\left(\ind_{\bw}(H_A(\bbeta))_{|\bv_0\to\bzero}\right)_{\bttheta\to\partial_{\by}})^{\perp,\C[\by]},
\end{equation}
and
\begin{equation}
\label{eq:indicial-colon-identification}
\left(\ind_{\bw}(H_A(\bbeta))_{|\bv_0\to\bzero}\right)_{\bttheta\to\partial_{\by}}=\Ann_{\C[\partial_{\by}]}^\bullet(C_{\bv_0,\cN_{\bv_0}})=\fkQ_{\cN_{\bv_0}}:\fkm_{\cN_{\bv_0}}.
\end{equation}
\end{theorem}
\begin{proof}
For brevity, put
\[
J:=\ind_{\bw}(H_A(\bbeta))\subset\C[\bttheta],\quad \fkJ_{\bv_0}:=\left(J_{|\bv_0\to\bzero}\right)_{\bttheta\to\partial_{\by}}
\subset\C[\partial_{\by}].
\]
Since $H_A(\bbeta)$ is holonomic and $\bw$ is generic, \cite[Theorem~2.3.9]{SST} shows that $J$ is a holonomic Frobenius ideal.  
In particular, it is zero-dimensional.
Hence its component supported at $\bv_0$, translated to the origin and identified with an ideal of $\C[\partial_{\by}]$, is Artinian.  
Thus $\fkJ_{\bv_0}$ is an Artinian ideal.

By \cite[Theorem~2.3.11]{SST}, the space of solutions of the component of $J$ supported at $\bv_0$ is precisely $\{\bx^{\bv_0}r(\log\bx)\,|\,r(\by)\in (\fkJ_{\bv_0})^{\perp,\C[\by]}\}$.

We prove \eqref{eq:Cv0-indicial-perp}.
First, let $r(\by)\in C_{\bv_0,\cN_{\bv_0}}$.
By the definition of $C_{\bv_0,\cN_{\bv_0}}$, there exists an
$A$-hypergeometric series in direction $\bw$ of the form
\[
\phi(\bx)=\sum_{\bu\in L, I_{\bu}\in \cN_{\bv_0}}\bx^{\bv_0+\bu}r_{\bu}(\log\bx)
\]
whose coefficient at $\bx^{\bv_0}$ is $r(\log\bx)$.
The definition of $\cN_{\bv_0}$ (see \eqref{eq:canonical-negative-support-family}) gives $\bw\cdot\bu\geq0$ for every $\bu\in L$ occurring in this series.
Since $\bw$ is generic, equality holds only for $\bu=\bzero$.
Consequently,
\[
\operatorname{in}_{\bw}(\phi)=\bx^{\bv_0}r(\log\bx).
\]
By \cite[Theorem~2.5.5]{SST}, $\operatorname{in}_{\bw}(\phi)$ is annihilated by $\operatorname{in}_{(-\bw,\bw)}(H_A(\bbeta))$.
The ideal $\operatorname{in}_{(-\bw,\bw)}(H_A(\bbeta))$ and $J$ have the same solutions among the logarithmic series under consideration; see the discussion preceding \cite[Lemma~2.5.10]{SST}.
It follows that $\bx^{\bv_0}r(\log\bx)$ is annihilated by the component of $J$ supported at $\bv_0$.
After translating $\bv_0$ to the origin and replacing
$\theta_j$ by $\partial_{y_j}$ for $j=1,\dots,n$, we obtain
\[
r(\by)\in(\fkJ_{\bv_0})^{\perp,\C[\by]}.
\]
Therefore,
\begin{equation}
\label{eq:Cv0-contained-indicial-perp-proof}
C_{\bv_0,\cN_{\bv_0}}\subset (\fkJ_{\bv_0})^{\perp,\C[\by]}.
\end{equation}

Conversely, let $r(\by)\in(\fkJ_{\bv_0})^{\perp,\C[\by]} =\Sol_{\C[\by]}^\bullet(\fkJ_{\bv_0})$. By \cite[Theorems~2.3.11 and~2.5.12]{SST}, the logarithmic term $\bx^{\bv_0}r(\log\bx)$ extends to an $A$-hypergeometric series in direction $\bw$. By \Cref{lem:direction-support-contained-canonical-family}, every negative support occurring in this series belongs to $\cN_{\bv_0}$. Hence $r(\by)\in C_{\bv_0,\cN_{\bv_0}}$, which proves the reverse inclusion. Together with \eqref{eq:Cv0-contained-indicial-perp-proof}, this proves \eqref{eq:Cv0-indicial-perp}.

Since $\fkJ_{\bv_0}$ is Artinian, \Cref{lem:artinian-ideal-orthogonal-duality} gives
\[
\Ann_{\C[\partial_{\by}]}^\bullet\left((\fkJ_{\bv_0})^{\perp,\C[\by]}\right)=\fkJ_{\bv_0}.
\]
Taking annihilators in \eqref{eq:Cv0-indicial-perp} and using
\eqref{eq:ann-CvN-colon}, we obtain
\[
\fkJ_{\bv_0}=\Ann_{\C[\partial_{\by}]}^\bullet\left((\fkJ_{\bv_0})^{\perp,\C[\by]}\right)=\Ann_{\C[\partial_{\by}]}^\bullet(C_{\bv_0,\cN_{\bv_0}})=\fkQ_{\cN_{\bv_0}}:\fkm_{\cN_{\bv_0}}.
\]
\end{proof}

\section{Ambient perturbation}
\label{sec:ambient-perturbation-general-N}

In this section we construct an ambient perturbation in the $\bt$-variables.
The $\bt$-variables belong to the formal series star-duality side of \Cref{subsec:t-variables-At-perp-comparison}.
Thus orthogonal complements of ideals in $\C[[\bt]]$ are taken in $\C[\partial_{\bt}]$, and we write them as $(-)^{\perp,\C[\partial_{\bt}]}$.

The goal of this section is to show that the ambient perturbation method realizes the coefficient space $C_{\bv_0,N}$ described in \Cref{thm:ann-CvN}.

\subsection{Ambient perturbation coefficients}
\label{subsec:ambient-perturbation-coefficients}

Let $\bt=(t_1,\dots,t_n)$.
For a subset $S\subset\{1,\dots,n\}$, we write
\[
\bt^S:=\prod_{j\in S}t_j,
\]
with the convention that the empty product is $1$.
We set
\begin{equation}
\label{eq:mhatN-definition}
\widehat m_N(\bt):=\bt^{I_{\bzero}\setminus K_N}\in\C[[\bt]].
\end{equation}

For $\bu\in L$, define
\begin{equation}
\label{eq:ambient-a-u}
\widehat a_{\bu}(\bt):=\frac{[\bv_0+\bt]_{\bu_-}}{[\bv_0+\bt+\bu]_{\bu_+}}.
\end{equation}
This is the ambient analogue of the intrinsic perturbation coefficient $a_\bu(\bs)$ as in \eqref{eq:intrinsic-coefficient-au}.
Here the exponent is perturbed in the ambient space $\C^n$, rather than inside the affine subspace $\bv_0+L_{\C}$.
The expression $\widehat a_{\bu}(\bt)$ is regarded as an element of the quotient field of $\C[[\bt]]$.
\Cref{lem:ambient-monomial-factor} below shows that, after multiplication by the normalizing monomial $\widehat m_N(\bt)$, the coefficient becomes a formal power series for every $\bu\in L$ with $I_{\bu}\in N$.

\begin{lemma}
\label{lem:ambient-coefficient-identity}
For every $\bu,\bu'\in L$, one has
\begin{equation}
\label{eq:ambient-transport}
\widehat a_{\bu+\bu'}(\bt)[\bv_0+\bt+\bu+\bu']_{\bu'_+}=\widehat a_{\bu}(\bt)[\bv_0+\bt+\bu]_{\bu'_-}.
\end{equation}
\end{lemma}

\begin{proof}
It is enough to prove the corresponding one-variable identity in each coordinate.
Fix $j\in\{1,\dots,n\}$, and put
\[
\alpha=(\bv_0)_j+t_j,\quad m=u_j,\quad n=u'_j.
\]
For $m\in\Z$, define
\[
R_m(\alpha):=\frac{[\alpha]_{m_-}}{[\alpha+m]_{m_+}}.
\]
Equivalently,
\begin{equation}
\label{eq:ambient-R-explicit}
R_m(\alpha)=
\begin{cases}
\displaystyle \prod_{k=1}^{m}(\alpha+k)^{-1}, & m\geq0,\\[8pt]
\displaystyle \prod_{k=m+1}^{0}(\alpha+k), & m<0,
\end{cases}
\end{equation}
where an empty product is understood to be $1$.
It follows directly from \eqref{eq:ambient-R-explicit}, after cancelling the common factors, that
\begin{equation}
\label{eq:ambient-R-product-explicit}
R_m(\alpha)R_n(\alpha+m)=
\begin{cases}
\displaystyle \prod_{k=1}^{m+n}(\alpha+k)^{-1}, & m+n\geq0,\\[8pt]
\displaystyle \prod_{k=m+n+1}^{0}(\alpha+k), & m+n<0.
\end{cases}
\end{equation}
Comparing \eqref{eq:ambient-R-product-explicit} with
\eqref{eq:ambient-R-explicit}, applied to $m+n$, gives
\begin{equation}
\label{eq:ambient-R-product}
R_{m+n}(\alpha)=R_m(\alpha)R_n(\alpha+m).
\end{equation}

By the definition of $R_n$, applied with $\alpha$ replaced by $\alpha+m$, we also have
\begin{equation}
\label{eq:ambient-R-shift}
R_n(\alpha+m)[\alpha+m+n]_{n_+}=[\alpha+m]_{n_-}.
\end{equation}
Multiplying \eqref{eq:ambient-R-shift} by $R_m(\alpha)$ and using \eqref{eq:ambient-R-product}, we get
\[
R_{m+n}(\alpha)[\alpha+m+n]_{n_+}=R_m(\alpha)[\alpha+m]_{n_-}.
\]
Taking the product of this identity over $j=1,\dots,n$ gives \eqref{eq:ambient-transport}.
\end{proof}

\begin{lemma}
\label{lem:ambient-monomial-factor}
For any $\bu\in L$ with $I_{\bu}\in N$, there exists a unit $\widehat b_{\bu}(\bt)\in\C[[\bt]]^\times$ such that
\begin{equation}
\label{eq:ambient-monomial-factor}
\widehat m_N(\bt)\widehat a_{\bu}(\bt)=\widehat b_{\bu}(\bt)\,\bt^{I_{\bu}\setminus K_N}.
\end{equation}
\end{lemma}

\begin{proof}
Fix $j\in\{1,\dots,n\}$.
The $j$-th factor of $\widehat a_{\bu}(\bt)$ is
\[
\frac{[(\bv_0)_j+t_j]_{(-u_j)_+}}{[(\bv_0)_j+t_j+u_j]_{(u_j)_+}}.
\]
A direct inspection of the consecutive affine factors gives
\begin{equation}
\label{eq:ambient-au-tj-order}
\ord_{t_j}\bigl(\widehat a_{\bu}(\bt)\bigr)=\mathbf 1_{I_{\bu}\setminus I_{\bzero}}(j)-\mathbf 1_{I_{\bzero}\setminus I_{\bu}}(j),
\end{equation}
where $\mathbf 1_S$ denotes the indicator function of a subset $S\subset\{1,\dots,n\}$.
Indeed, a factor $t_j$ occurs in the numerator precisely when $j\in I_{\bu}\setminus I_{\bzero}$, and a factor $t_j$ occurs in the denominator precisely when $j\in I_{\bzero}\setminus I_{\bu}$.
In either case its multiplicity is one.

Since $\widehat m_N(\bt)=\bt^{I_{\bzero}\setminus K_N}$ and $K_N\subset I_{\bzero}\cap I_{\bu}$, it follows from
\eqref{eq:ambient-au-tj-order} that
\begin{align}
\ord_{t_j}\bigl(\widehat m_N(\bt)\widehat a_{\bu}(\bt)\bigr)
&=\mathbf 1_{I_{\bzero}\setminus K_N}(j)+\mathbf 1_{I_{\bu}\setminus I_{\bzero}}(j)-\mathbf 1_{I_{\bzero}\setminus I_{\bu}}(j)\notag\\
&=\mathbf 1_{I_{\bu}\setminus K_N}(j).
\label{eq:ambient-normalized-tj-order}
\end{align}
Thus $\widehat m_N(\bt)\widehat a_{\bu}(\bt)$ has exactly the monomial factor $\bt^{I_{\bu}\setminus K_N}$. 
After removing this factor, the remaining formal power series has nonzero constant term and is therefore a unit in $\C[[\bt]]$. 
Hence the assertion holds.
\end{proof}

\subsection{The ambient perturbation series}
\label{subsec:ambient-perturbation-series}

The \emph{ambient perturbation series associated with $N$} is
\begin{equation}
\label{eq:ambient-perturbation-series}
\widehat F_N(\bx,\bt):=\widehat m_N(\bt)\sum_{\bu\in L,\ I_{\bu}\in N}\widehat a_{\bu}(\bt)\,\bx^{\bv_0+\bt+\bu}.
\end{equation}
By \Cref{lem:ambient-monomial-factor}, this is a well-defined formal series in the $\bt$-variables.
More precisely, it may be written as
\begin{equation}
\label{eq:ambient-perturbation-expanded}
\widehat F_N(\bx,\bt)=\sum_{\bu\in L,\ I_{\bu}\in N}\widehat b_{\bu}(\bt)\,\bt^{I_{\bu}\setminus K_N}\bx^{\bv_0+\bt+\bu},\quad\widehat b_{\bu}(\bt)\in\C[[\bt]]^\times.
\end{equation}

Define the polynomial ideals
\[
\widehat M_N^{\mathrm{pol}}:=\left\langle\bt^{I\setminus K_N}\,\middle|\,I\in N\right\rangle\subset\C[\bt],
\]
and
\[
\widehat P_N^{\mathrm{pol}}:=\left\langle\bt^{(I\cup J)\setminus K_N}\,\middle|\,I\in N,\ J\in N^c\right\rangle\subset\C[\bt].
\]
Put
\begin{equation}
\label{eq:QhatN-pol-definition}
\widehat Q_N^{\mathrm{pol}}:=\langle A\bt\rangle_{\mathrm{pol}}\widehat M_N^{\mathrm{pol}}+\widehat P_N^{\mathrm{pol}}\subset\C[\bt],
\end{equation}
where
\[
\langle A\bt\rangle_{\mathrm{pol}}:=\left\langle\ba^{(i)}\bt\,\middle|\,i=1,\dots,d\right\rangle\subset\C[\bt].
\]
Its extension to $\C[[\bt]]$ is the ideal $\langle A\bt\rangle$ introduced in \eqref{eq:At}.

We denote their extensions to $\C[[\bt]]$ by
\[
\widehat M_N:=\widehat M_N^{\mathrm{pol}}\C[[\bt]],\quad\widehat P_N:=\widehat P_N^{\mathrm{pol}}\C[[\bt]],
\]
and
\begin{equation}
\label{eq:QhatN-formal-definition}
\widehat Q_N:=\widehat Q_N^{\mathrm{pol}}\C[[\bt]]=\langle A\bt\rangle\widehat M_N+\widehat P_N\subset\C[[\bt]].
\end{equation}
Then $\widehat Q_N$ is homogeneous.
Unless a coefficient ring is explicitly indicated, a colon ideal involving $\widehat Q_N$ is formed in $\C[[\bt]]$, whereas a colon ideal involving $\widehat Q_N^{\mathrm{pol}}$ is formed in $\C[\bt]$.

\begin{lemma}
\label{lem:ambient-H-action}
The following hold.
\begin{enumerate}[(1)]
\item
For each Euler operator $(A\theta-\bbeta)_i$ ($i=1,\dots,d$), every $\bt$-coefficient of $(A\theta-\bbeta)_i\bullet \widehat F_N(\bx,\bt)$ belongs to $\widehat Q_N$.

\item For each $\ggamma\in L$, every $\bt$-coefficient of $\left(\partial_{\bx}^{\ggamma_+}-\partial_{\bx}^{\ggamma_-}\right)\bullet \widehat F_N(\bx,\bt)$ belongs to $\widehat Q_N$.
\end{enumerate}
\end{lemma}

\begin{proof}
(1) For $i=1,\dots,d$, using $A(\bv_0+\bu)=\bbeta$ for all $\bu\in L$, we have
\[
(A\theta-\bbeta)_i\bullet \bx^{\bv_0+\bt+\bu}=(A\bt)_i\,\bx^{\bv_0+\bt+\bu}.
\]
Hence, by \eqref{eq:ambient-perturbation-expanded}, every coefficient appearing in $(A\theta-\bbeta)_i\bullet\widehat F_N(\bx,\bt)$ belongs to
\[
\langle A\bt\rangle\widehat M_N\subset\widehat Q_N.
\]

(2) Fix $\ggamma\in L$, and consider the toric operator $\partial_{\bx}^{\ggamma_+}-\partial_{\bx}^{\ggamma_-}$.
For $\bu\in L$, the two possible contributions to the monomial
\[
\bx^{\bv_0+\bt+\bu-\ggamma_+}=\bx^{\bv_0+\bt+(\bu-\ggamma)-\ggamma_-}
\]
have coefficients
\begin{equation}
\label{eq:ambient-plus-coefficient}
\widehat m_N(\bt)\widehat a_{\bu}(\bt)[\bv_0+\bt+\bu]_{\ggamma_+}
\end{equation}
and
\begin{equation}
\label{eq:ambient-minus-coefficient}
\widehat m_N(\bt)\widehat a_{\bu-\ggamma}(\bt)[\bv_0+\bt+\bu-\ggamma]_{\ggamma_-},
\end{equation}
respectively, whenever the corresponding negative supports belong to $N$.

If both $I_{\bu}$ and $I_{\bu-\ggamma}$ belong to $N$, then the two
coefficients are equal by \Cref{lem:ambient-coefficient-identity}, applied to $\bu-\ggamma$ and $\ggamma$.
Hence they cancel.

Suppose next that exactly one of $I_{\bu}$ and
$I_{\bu-\ggamma}$ belongs to $N$.
Assume that $I_{\bu}\in N$ and $I_{\bu-\ggamma}\in N^c$.
We show that the surviving coefficient \eqref{eq:ambient-plus-coefficient} belongs to $\widehat P_N$.
The case in which $I_{\bu}\in N^c$ and $I_{\bu-\ggamma}\in N$ is proved in the same way, with the two indices interchanged.
By \Cref{lem:ambient-monomial-factor}, $\widehat m_N(\bt)\widehat a_{\bu}(\bt)$ is a unit times $\bt^{I_{\bu}\setminus K_N}$.
If $j\in I_{\bu-\ggamma}\setminus I_{\bu}$, then $(\bv_0+\bu)_j\in\N$ and $(\bv_0+\bu-\ggamma)_j\in\Z_{<0}$.
Hence $\gamma_j>0$, and the $j$-th factor of $[\bv_0+\bt+\bu]_{\ggamma_+}$ is divisible by $t_j$.
Therefore, \eqref{eq:ambient-plus-coefficient} is divisible by
\[
\bt^{I_{\bu}\setminus K_N}\bt^{I_{\bu-\ggamma}\setminus I_{\bu}}=\bt^{(I_{\bu}\cup I_{\bu-\ggamma})\setminus K_N}.
\]
Since $I_{\bu}\in N$ and $I_{\bu-\ggamma}\in N^c$, the monomial on the right-hand side is a generator of $\widehat P_N$.
Thus the surviving coefficient belongs to $\widehat P_N\subset\widehat Q_N$.

If both $I_{\bu}$ and $I_{\bu-\ggamma}$ belong to $N^c$, neither summand occurs.
Consequently, every $\bt$-coefficient of $\left(\partial_{\bx}^{\ggamma_+}-\partial_{\bx}^{\ggamma_-}\right)\bullet\widehat F_N(\bx,\bt)$ belongs to $\widehat Q_N$.
\end{proof}

\subsection{The ambient coefficient space}
\label{subsec:ambient-coefficient-space}

We now identify the coefficient space obtained from ambient perturbation.
Define the monomial-preserving $\C$-linear isomorphism
\[
\Xi:\C[\partial_{\bt}]\longrightarrow \C[\by],\quad\Xi(\partial_{\bt}^{\aalpha})=\by^{\aalpha}.
\]
We also define a $\C$-algebra isomorphism
\[
\Theta:\C[\bt]\longrightarrow \C[\partial_{\by}],\quad\Theta(\bt^{\aalpha})=\partial_{\by}^{\aalpha}.
\]
Then
\begin{equation}
\label{eq:Theta-Qhat-Qfrak}
\Theta(\widehat m_N(\bt))=\fkm_N,\ \ \Theta(\widehat M_N^{\mathrm{pol}})=\fkM_N,\ \ \Theta(\widehat P_N^{\mathrm{pol}})=\fkP_N,\ \ \Theta(\widehat Q_N^{\mathrm{pol}})=\fkQ_N.
\end{equation}
Since $\Theta$ is an algebra isomorphism, it commutes with colon ideals.
Hence
\begin{equation}
\label{eq:Theta-colon}
\Theta\left(\widehat Q_N^{\mathrm{pol}}:\widehat m_N(\bt)\right)=\fkQ_N:\fkm_N.
\end{equation}

Moreover, flatness of the completion
$\C[\bt]\to\C[[\bt]]$ gives
\begin{equation}
\label{eq:colon-commutes-with-completion}
\left(\widehat Q_N^{\mathrm{pol}}:\widehat m_N(\bt)\right)\C[[\bt]]=\widehat Q_N:\widehat m_N(\bt).
\end{equation}
Thus the formal series colon ideal appearing in the star-duality construction is the extension of the polynomial colon ideal occurring in \eqref{eq:Theta-colon}.

\begin{lemma}
\label{lem:Xi-perp-Theta}
Let $J\subset\C[\bt]$ be a homogeneous ideal, and let $\C[[\bt]]J$ denote its extension to $\C[[\bt]]$.
Then
\begin{equation}
\label{eq:Xi-perp-Theta}
\Xi\left((\C[[\bt]]J)^{\perp,\C[\partial_{\bt}]}\right)=\Theta(J)^{\perp,\C[\by]}.
\end{equation}
\end{lemma}

\begin{proof}
For $\aalpha,\bbeta\in\N^n$, we have
\begin{equation}
\label{eq:Xi-Theta-monomial-pairing}
(\partial_{\bt}^{\aalpha},\bt^{\bbeta})_{\bt}=(\partial_{\by}^{\bbeta},\by^{\aalpha})_{\by}=\delta_{\aalpha,\bbeta}\,\aalpha!.
\end{equation}
Hence, for every $q\in\C[\partial_{\bt}]$ and every $f\in\C[\bt]$,
\begin{equation}
\label{eq:Xi-Theta-pairing-compatibility}
(q,f)_{\bt}=(\Theta(f),\Xi(q))_{\by}.
\end{equation}

For $q\in\C[\partial_{\bt}]$, orthogonality to $\C[[\bt]]J$ is equivalent to
\[
(q,\bt^{\ggamma}h)_{\bt}=0\quad(\ggamma\in\N^n,\ h\in J).
\]
By \eqref{eq:Xi-Theta-pairing-compatibility}, this is equivalent to
\[
(\partial_{\by}^{\ggamma},\Theta(h)\bullet\Xi(q))_{\by}=0\quad(\ggamma\in\N^n,\ h\in J).
\]
By non-degeneracy of the ordinary pairing, the latter condition holds if and only if
\[
\Theta(h)\bullet\Xi(q)=0\quad(h\in J).
\]
This proves the assertion.
\end{proof}

\begin{theorem}
\label{thm:ambient-perturbation-general-N}
Let $N$ be an ordered negative support family for $\bv_0$.
For any $q(\partial_{\bt})\in(\widehat Q_N)^{\perp,\C[\partial_{\bt}]}$, the series
\begin{equation}
\label{eq:ambient-series-output}
\left(q(\partial_{\bt})\bullet \widehat F_N(\bx,\bt)\right)_{|\bt=\bzero}
\end{equation}
is an $A$-hypergeometric series in direction $\bw$.

Moreover, the coefficients of $\bx^{\bv_0}$ obtained in this way are precisely
\begin{equation}
\label{eq:ambient-coefficient-space}
\Xi\left((\widehat Q_N:\widehat m_N(\bt))^{\perp,\C[\partial_{\bt}]}\right)
\end{equation}
and
\begin{equation}
\label{eq:ambient-realizes-CvN}
C_{\bv_0,N}=\Xi\left((\widehat Q_N:\widehat m_N(\bt))^{\perp,\C[\partial_{\bt}]}\right).
\end{equation}
\end{theorem}

\begin{proof}
Let $q\in\widehat Q_N^{\perp,\C[\partial_{\bt}]}$.
By \Cref{lem:Pperp}, this is equivalent to $f\star q=0$ for all $f\in\widehat Q_N$.
By \Cref{lem:ambient-H-action}, every $\bt$-coefficient of the action of an Euler or toric generator of $H_A(\bbeta)$ on $\widehat F_N(\bx,\bt)$ belongs to $\widehat Q_N$.
Applying $q(\partial_{\bt})$ and evaluating at $\bt=\bzero$ therefore annihilate all these coefficients.
Hence $\left(q(\partial_{\bt})\bullet\widehat F_N(\bx,\bt)\right)_{|\bt=\bzero}$ is annihilated by $H_A(\bbeta)$.

We next compute the coefficient of $\bx^{\bv_0}$.
Among the summands in the defining series \eqref{eq:ambient-perturbation-series}, only the summand indexed by $\bu=\bzero$ can contribute to this coefficient after applying $q(\partial_{\bt})$ and evaluating at $\bt=\bzero$.
Moreover, $\widehat a_{\bzero}(\bt)=1$.
By the same calculation as in
\eqref{eq:intrinsic-leading-coefficient-formula}, we have
\begin{equation}
\label{eq:ambient-leading-coefficient-formula}
\left(q(\partial_{\bt})\bullet\bigl(\widehat m_N(\bt)\bx^{\bv_0+\bt}\bigr)\right)_{|\bt=\bzero}=\bx^{\bv_0}\Xi\bigl(\widehat m_N(\bt)\star q(\partial_{\bt})\bigr)(\log\bx).
\end{equation}
Therefore, the coefficient of $\bx^{\bv_0}$ obtained from $q$ is
\[
\Xi\bigl(\widehat m_N(\bt)\star q(\partial_{\bt})\bigr)(\log\bx).
\]
Consequently, the full ambient coefficient space is
\[
\Xi\left(
\widehat m_N(\bt)\star
(\widehat Q_N)^{\perp,\C[\partial_{\bt}]}
\right).
\]
Since $\widehat Q_N$ and $\widehat m_N(\bt)$ are homogeneous, \Cref{prop:ColonIdeal} gives
\[
\widehat m_N(\bt)\star (\widehat Q_N)^{\perp,\C[\partial_{\bt}]}=(\widehat Q_N:\widehat m_N(\bt))^{\perp,\C[\partial_{\bt}]}.
\]
This proves \eqref{eq:ambient-coefficient-space}.

Finally, by \eqref{eq:colon-commutes-with-completion}, the formal series ideal $\widehat Q_N:\widehat m_N(\bt)$ is the extension of $\widehat Q_N^{\mathrm{pol}}:\widehat m_N(\bt)$.
Hence \Cref{lem:Xi-perp-Theta} and \eqref{eq:Theta-colon} give
\[
\Xi\left((\widehat Q_N:\widehat m_N(\bt))^{\perp,\C[\partial_{\bt}]}\right)=(\fkQ_N:\fkm_N)^{\perp,\C[\by]}.
\]
By \Cref{thm:ann-CvN}, the right-hand side is $C_{\bv_0,N}$.
Hence \eqref{eq:ambient-realizes-CvN} follows.
\end{proof}

\section{Comparison of intrinsic and ambient perturbations}
\label{sec:comparison-intrinsic-ambient-perturbations}

In this section, we compare the coefficient spaces obtained by the intrinsic and ambient perturbation methods.
We keep the distinction between the $\bs$-side and the $\bt$-side notation.

\subsection{The comparison criterion}
\label{subsec:comparison-criterion}

By \Cref{thm:ambient-perturbation-general-N}, the ambient perturbation method realizes the coefficient space
\begin{equation}
\label{eq:ambient-space-section6}
V_N^{\mathrm{amb}}:=\Xi\left((\widehat Q_N:\widehat m_N(\bt))^{\perp,\C[\partial_{\bt}]}\right)=C_{\bv_0,N}.
\end{equation}

On the other hand, by \Cref{thm:intrinsic-perturbation-general-N}, the coefficients of $\bx^{\bv_0}$ obtained by intrinsic perturbation form
\begin{equation}
\label{eq:intrinsic-space-section6}
V_N^{\mathrm{int}}:=\Xi\Psi\left((P_N:m_N(\bs))^{\perp,\C[\partial_{\bs}]}\right)\subset\C[\by].
\end{equation}
Here $\Psi:\C[\partial_{\bs}]\to\C[\partial_{\bt}]$ is the map introduced in \Cref{subsec:t-variables-At-perp-comparison}, and $\Xi:\C[\partial_{\bt}]\to\C[\by]$ is the map introduced in \Cref{subsec:ambient-coefficient-space}.

\begin{lemma}
\label{lem:Jhat-contains-At-and-inclusion}
The ideal $\widehat Q_N:\widehat m_N(\bt)$ contains $\langle A\bt\rangle$, and
\begin{equation}
\label{eq:Phi-Jhat-in-J}
\Phi(\widehat Q_N:\widehat m_N(\bt))\subset P_N:m_N(\bs).
\end{equation}
\end{lemma}

\begin{proof}
Since $\widehat m_N(\bt)\in\widehat M_N$, we have $\langle A\bt\rangle\,\widehat m_N(\bt)\subset\langle A\bt\rangle\widehat M_N\subset\widehat Q_N$.
Hence $\langle A\bt\rangle\subset \widehat Q_N:\widehat m_N(\bt)$.

Let $f\in\widehat Q_N:\widehat m_N(\bt)$.
Then, using \eqref{eq:QhatN-formal-definition},
\[
f\,\widehat m_N(\bt)\in
\widehat Q_N=\langle A\bt\rangle\widehat M_N+\widehat P_N.
\]
Applying $\Phi$, and using $\Phi(\langle A\bt\rangle)=0$, $\Phi(\widehat m_N(\bt))=m_N(\bs)$, and $\Phi(\widehat P_N)=P_N$, we obtain $\Phi(f)m_N(\bs)\in P_N$.
Thus $\Phi(f)\in P_N:m_N(\bs)$.
This proves \eqref{eq:Phi-Jhat-in-J}.
\end{proof}

\begin{proposition}
\label{prop:L-perturb-general-N}
The intrinsic coefficient space is contained in the ambient coefficient space:
\begin{equation}
\label{eq:int-contained-amb}
V_N^{\mathrm{int}}\subset V_N^{\mathrm{amb}}.
\end{equation}
Moreover, the following conditions are equivalent.
\begin{enumerate}[(1)]
\item The intrinsic perturbation method realizes the full coefficient space, that is,
\begin{equation}
\label{eq:int-equals-amb-space}
V_N^{\mathrm{int}}=V_N^{\mathrm{amb}}.
\end{equation}

\item The colon ideals satisfy
\begin{equation}
\label{eq:Phi-Jhat-equals-J}
\Phi(\widehat Q_N:\widehat m_N(\bt))=P_N:m_N(\bs).
\end{equation}
\end{enumerate}
\end{proposition}

\begin{proof}
By \Cref{lem:Jhat-contains-At-and-inclusion}, the ideal $\widehat Q_N:\widehat m_N(\bt)$ contains $\langle A\bt\rangle$.
Hence \Cref{prop:Ay-perp2} applies to $\widehat Q_N:\widehat m_N(\bt)$, and gives
\begin{equation}
\label{eq:Jhat-perp-via-Phi}
(\widehat Q_N:\widehat m_N(\bt))^{\perp,\C[\partial_{\bt}]}=\Psi\left(\Phi(\widehat Q_N:\widehat m_N(\bt))^{\perp,\C[\partial_{\bs}]}\right).
\end{equation}
Therefore, by \eqref{eq:ambient-space-section6},
\begin{equation}
\label{eq:amb-space-via-Phi-Jhat}
V_N^{\mathrm{amb}}=\Xi\Psi\left(\Phi(\widehat Q_N:\widehat m_N(\bt))^{\perp,\C[\partial_{\bs}]}\right).
\end{equation}
On the other hand, by \eqref{eq:intrinsic-space-section6},
\[
V_N^{\mathrm{int}}=\Xi\Psi\left((P_N:m_N(\bs))^{\perp,\C[\partial_{\bs}]}\right).
\]
Since $\Phi(\widehat Q_N:\widehat m_N(\bt))\subset P_N:m_N(\bs)$ by \Cref{lem:Jhat-contains-At-and-inclusion}, we have
\[
(P_N:m_N(\bs))^{\perp,\C[\partial_{\bs}]}\subset\Phi(\widehat Q_N:\widehat m_N(\bt))^{\perp,\C[\partial_{\bs}]}.
\]
Applying the injective map $\Xi\Psi$, we obtain $V_N^{\mathrm{int}}\subset V_N^{\mathrm{amb}}$.

We prove the equivalence.
If $\Phi(\widehat Q_N:\widehat m_N(\bt))=
P_N:m_N(\bs)$, then \eqref{eq:amb-space-via-Phi-Jhat} immediately gives $V_N^{\mathrm{amb}}=V_N^{\mathrm{int}}$.

Conversely, suppose that $V_N^{\mathrm{amb}}=V_N^{\mathrm{int}}$.
Since $\Xi\Psi$ is injective, it follows from \eqref{eq:amb-space-via-Phi-Jhat} and \eqref{eq:intrinsic-space-section6} that
\[
\Phi(\widehat Q_N:\widehat m_N(\bt))^{\perp,\C[\partial_{\bs}]}=(P_N:m_N(\bs))^{\perp,\C[\partial_{\bs}]}.
\]
Both $\Phi(\widehat Q_N:\widehat m_N(\bt))$ and $P_N:m_N(\bs)$ are homogeneous ideals of $\C[[\bs]]$.
Taking double orthogonal complements and using \Cref{prop:Duality}, we obtain
\[
\Phi(\widehat Q_N:\widehat m_N(\bt))=P_N:m_N(\bs).
\]
\end{proof}

In the rest of this section, once the ideal equality \eqref{eq:Phi-Jhat-equals-J} is established, the equality
\[
V_N^{\mathrm{int}}=V_N^{\mathrm{amb}}=C_{\bv_0,N}
\]
follows immediately from \Cref{prop:L-perturb-general-N} and \eqref{eq:ambient-space-section6}.  We shall use this implication without further comment.

\subsection{A first sufficient condition}
\label{subsec:first-sufficient-condition}

We first record a simple sufficient condition for the comparison equality \eqref{eq:Phi-Jhat-equals-J}.
It is expressed in terms of the ideal determined by the fixed support $I_{\bzero}$.

Define
\begin{equation}
\label{eq:PN-I0-definition}
P_{N,I_{\bzero}}:=\left\langle(B\bs)^{J\setminus I_{\bzero}}\,\middle|\,J\in N^c\right\rangle\subset\C[[\bs]],
\end{equation}
and define its ambient counterpart by
\begin{equation}
\label{eq:PhatN-I0-definition}
\widehat P_{N,I_{\bzero}}:=\left\langle\bt^{J\setminus I_{\bzero}}\,\middle|\,J\in N^c\right\rangle\subset\C[[\bt]].
\end{equation}
Since $\Phi$ is surjective, one has
\begin{equation}
\label{eq:Phi-PhatNI0}
\Phi(\widehat P_{N,I_{\bzero}})=P_{N,I_{\bzero}}.
\end{equation}

\begin{lemma}
\label{lem:PN-colon-contained-pI0}
The following inclusions hold:
\begin{equation}
\label{eq:Phat-inclusions-I0}
\langle A\bt\rangle+\widehat P_{N,I_{\bzero}}\subset\widehat Q_N:\widehat m_N(\bt),
\end{equation}
and
\begin{equation}
\label{eq:P-inclusions-I0}
P_{N,I_{\bzero}}\subset P_N:m_N(\bs).
\end{equation}
\end{lemma}

\begin{proof}
We first prove \eqref{eq:Phat-inclusions-I0}.
As we have proved $\langle A\bt\rangle\subset\widehat Q_N:\widehat m_N(\bt)$ in \Cref{lem:Jhat-contains-At-and-inclusion}, it remains to show
\[
\widehat P_{N,I_{\bzero}}\subset\widehat Q_N:\widehat m_N(\bt).
\]
Let $J\in N^c$.
Since $\widehat m_N(\bt)=\bt^{I_{\bzero}\setminus K_N}$, we have
\[
\bt^{J\setminus I_{\bzero}}\widehat m_N(\bt)=\bt^{J\setminus I_{\bzero}}\bt^{I_{\bzero}\setminus K_N}=\bt^{(I_{\bzero}\cup J)\setminus K_N}.
\]
Since $I_{\bzero}\in N$ and $J\in N^c$, the last monomial is one of the generators of $\widehat P_N$.
Hence it belongs to $\widehat Q_N$, and therefore $\bt^{J\setminus I_{\bzero}}\in\widehat Q_N:\widehat m_N(\bt)$.
This proves \eqref{eq:Phat-inclusions-I0}.

The proof of \eqref{eq:P-inclusions-I0} is the $\bs$-side analogue.
For $J\in N^c$, we have
\[
(B\bs)^{J\setminus I_{\bzero}}\,m_N(\bs)=(B\bs)^{J\setminus I_{\bzero}}(B\bs)^{I_{\bzero}\setminus K_N}=(B\bs)^{(I_{\bzero}\cup J)\setminus K_N}.
\]
By the definition \eqref{eq:PNs} of $P_N$, this element is a generator of $P_N$.
Hence $(B\bs)^{J\setminus I_{\bzero}}\in P_N:m_N(\bs)$.
Thus $P_{N,I_{\bzero}}\subset P_N:m_N(\bs)$.
\end{proof}

\begin{proposition}
\label{prop:colon-y-general-N}
Assume that
\begin{equation}
\label{eq:first-sufficient-assumption}
P_N:m_N(\bs)=P_{N,I_{\bzero}}.
\end{equation}
Then
\begin{equation}
\label{eq:first-sufficient-comparison}
\Phi(\widehat Q_N:\widehat m_N(\bt))=P_N:m_N(\bs).
\end{equation}
\end{proposition}

\begin{proof}
By \eqref{eq:Phat-inclusions-I0}, we have $\widehat P_{N,I_{\bzero}}\subset\widehat Q_N:\widehat m_N(\bt)$.
Applying $\Phi$ and using \eqref{eq:Phi-PhatNI0}, we obtain
\begin{equation}
\label{eq:equation-colon-sum}
P_{N,I_{\bzero}}=\Phi(\widehat P_{N,I_{\bzero}})\subset\Phi(\widehat Q_N:\widehat m_N(\bt)).
\end{equation}
On the other hand, by \Cref{lem:Jhat-contains-At-and-inclusion}, $\Phi(\widehat Q_N:\widehat m_N(\bt))\subset P_N:m_N(\bs)$.
Using the assumption \eqref{eq:first-sufficient-assumption}, we get $\Phi(\widehat Q_N:\widehat m_N(\bt))\subset P_{N,I_{\bzero}}$.
Together with \eqref{eq:equation-colon-sum}, this gives
\[
\Phi(\widehat Q_N:\widehat m_N(\bt))=P_{N,I_{\bzero}}=P_N:m_N(\bs),
\]
which proves \eqref{eq:first-sufficient-comparison}, equivalently \eqref{eq:Phi-Jhat-equals-J} under the assumption \eqref{eq:first-sufficient-assumption}.
\end{proof}

\subsection{Equality criteria}
\label{subsec:equality-criteria}

We now give more flexible criteria for the comparison equality \eqref{eq:Phi-Jhat-equals-J}.
The key point is to compare $\widehat Q_N:\widehat m_N(\bt)$ with the pullback of $P_N:m_N(\bs)$ under $\Phi$.

In this subsection, whenever an ideal generated in $\C[\bt]$ is used with $\Phi^{-1}$, intersections, or colon operations inside $\C[[\bt]]$, we regard it as extended to $\C[[\bt]]$.

\begin{lemma}
\label{lem:pullback-JN}
The pullback of $P_N:m_N(\bs)$ by $\Phi$ is
\begin{equation}
\label{eq:pullback-JN}
\Phi^{-1}(P_N:m_N(\bs))=\bigl(\langle A\bt\rangle\cap\widehat M_N+\widehat P_N\bigr):\widehat m_N(\bt).
\end{equation}
\end{lemma}

\begin{proof}
Let $f\in\C[[\bt]]$.
We have
\[
f\in\Phi^{-1}(P_N:m_N(\bs))\iff\Phi(f)m_N(\bs)\in P_N.
\]
Since $m_N(\bs)=\Phi(\widehat m_N(\bt))$, we see that
\[
\Phi\bigl(f\,\widehat m_N(\bt)\bigr)\in P_N\iff f\,\widehat m_N(\bt)\in \Phi^{-1}(P_N).
\]

We now compute $\Phi^{-1}(P_N)$.
Since $P_N=\Phi(\widehat P_N)$ and $\Ker(\Phi)=\langle A\bt\rangle$, we have $\Phi^{-1}(P_N)=\langle A\bt\rangle+\widehat P_N$.
Therefore, 
\[
f\in\Phi^{-1}(P_N:m_N(\bs))\iff f\,\widehat m_N(\bt)\in\langle A\bt\rangle+\widehat P_N.
\]

On the other hand, $\widehat m_N(\bt)\in\widehat M_N$ and $\widehat M_N$ is an ideal.
Hence $f\,\widehat m_N(\bt)\in\widehat M_N$.
Since $\widehat P_N\subset\widehat M_N$, we get
\[
\bigl(\langle A\bt\rangle+\widehat P_N\bigr)\cap\widehat M_N=\langle A\bt\rangle\cap\widehat M_N+\widehat P_N.
\]
Thus the preceding condition is equivalent to
\[
f\,\widehat m_N(\bt)\in \langle A\bt\rangle\cap\widehat M_N+\widehat P_N,
\]
which is precisely
\[
f\in\bigl(\langle A\bt\rangle\cap\widehat M_N+\widehat P_N\bigr):\widehat m_N(\bt).
\]
This proves \eqref{eq:pullback-JN}.
\end{proof}

\begin{proposition}
\label{prop:intersection-visible-general-N}
Assume that
\begin{equation}
\label{eq:intersection-visible-assumption}
\bigl(\langle A\bt\rangle\widehat M_N+\widehat P_N\bigr):\widehat m_N(\bt)=\bigl(\langle A\bt\rangle\cap\widehat M_N+\widehat P_N\bigr):\widehat m_N(\bt).
\end{equation}
Then
\begin{equation}
\label{eq:intersection-visible-conclusion}
\Phi(\widehat Q_N:\widehat m_N(\bt))=P_N:m_N(\bs).
\end{equation}
\end{proposition}

\begin{proof}
By \eqref{eq:QhatN-formal-definition}, the left-hand side of \eqref{eq:intersection-visible-assumption} is $\widehat Q_N:\widehat m_N(\bt)$.
By \Cref{lem:pullback-JN}, the right-hand side of \eqref{eq:intersection-visible-assumption} is $\Phi^{-1}(P_N:m_N(\bs))$.
Therefore the assumption gives
\[
\widehat Q_N:\widehat m_N(\bt)=\Phi^{-1}(P_N:m_N(\bs)).
\]
Applying $\Phi$ to both sides, and using the surjectivity of $\Phi$, we obtain
\[
\Phi(\widehat Q_N:\widehat m_N(\bt))=P_N:m_N(\bs).
\]
This proves \eqref{eq:intersection-visible-conclusion}, hence the comparison equality \eqref{eq:Phi-Jhat-equals-J}.
\end{proof}

\begin{corollary}
\label{cor:product-intersection-general-N}
If
\begin{equation}
\label{eq:product-intersection-condition}
\langle A\bt\rangle\cap\widehat M_N=\langle A\bt\rangle\widehat M_N,
\end{equation}
then
\begin{equation}
\label{eq:product-intersection-conclusion}
\Phi(\widehat Q_N:\widehat m_N(\bt))=P_N:m_N(\bs).
\end{equation}
\end{corollary}

\begin{proof}
The condition \eqref{eq:product-intersection-condition} implies
\[
\langle A\bt\rangle\cap\widehat M_N+\widehat P_N=\langle A\bt\rangle\widehat M_N+\widehat P_N.
\]
Therefore the two colon ideals in \eqref{eq:intersection-visible-assumption} are equal.
The assertion follows from
\Cref{prop:intersection-visible-general-N}.
\end{proof}

\begin{proposition}
\label{prop:square-general-N}
Suppose that there exists $I\in N$ such that
\begin{equation}
\label{eq:square-condition-general-N}
\widehat Q_N:\bigl(\widehat m_N(\bt)\,\bt^{I\setminus K_N}\bigr)=\widehat Q_N:\widehat m_N(\bt).
\end{equation}
Then
\begin{equation}
\label{eq:square-condition-conclusion}
\Phi(\widehat Q_N:\widehat m_N(\bt))=P_N:m_N(\bs).
\end{equation}
\end{proposition}

\begin{proof}
By \Cref{prop:intersection-visible-general-N}, it is enough to prove
\begin{equation}
\label{eq:square-proof-target}
\bigl(\langle A\bt\rangle\cap\widehat M_N+\widehat P_N\bigr):\widehat m_N(\bt)\subset\widehat Q_N:\widehat m_N(\bt).
\end{equation}
The reverse inclusion is automatic because
\[
\widehat Q_N=\langle A\bt\rangle\widehat M_N+\widehat P_N\subset\langle A\bt\rangle\cap\widehat M_N+\widehat P_N.
\]

Let $f\in\bigl(\langle A\bt\rangle\cap\widehat M_N+\widehat P_N\bigr):\widehat m_N(\bt)$.
Then $f\,\widehat m_N(\bt)\in\langle A\bt\rangle\cap\widehat M_N+\widehat P_N$.
Multiplying by $\bt^{I\setminus K_N}$, we obtain
\[
f\,\widehat m_N(\bt)\,\bt^{I\setminus K_N}\in\bt^{I\setminus K_N}\bigl(\langle A\bt\rangle\cap\widehat M_N+\widehat P_N\bigr).
\]
Since $I\in N$, the monomial $\bt^{I\setminus K_N}$ is one of the generators of $\widehat M_N$.
Hence $\bt^{I\setminus K_N}(\langle A\bt\rangle\cap\widehat M_N)\subset\langle A\bt\rangle\widehat M_N$.
Moreover, $\bt^{I\setminus K_N}\widehat P_N\subset\widehat P_N$ because $\widehat P_N$ is an ideal.
Therefore
\[
f\,\widehat m_N(\bt)\,\bt^{I\setminus K_N}\in\langle A\bt\rangle\widehat M_N+\widehat P_N=\widehat Q_N.
\]
Thus $f\in\widehat Q_N:\bigl(\widehat m_N(\bt)\,\bt^{I\setminus K_N}\bigr)$.
By the assumption \eqref{eq:square-condition-general-N}, this gives $f\in\widehat Q_N:\widehat m_N(\bt)$.
This proves \eqref{eq:square-proof-target}. 
\end{proof}

\section{Examples}
\label{sec:examples}

We give two examples illustrating the comparison criterion of
\Cref{prop:L-perturb-general-N}.
In both examples the equality
\begin{equation}
\label{eq:examples-comparison-equality}
\Phi(\widehat Q_N:\widehat m_N(\bt))=P_N:m_N(\bs)
\end{equation}
holds.
Hence the intrinsic perturbation method realizes the full coefficient space.
The negative support data used below are recorded in \Cref{app:negative-support-computations}.

\begin{example}[{\cite[Example 3.14]{SST}}]
\label{ex:2.1}
Let
\[
A=
\begin{pmatrix}
  1&1&1&1\\
  0&1&2&3
\end{pmatrix},
\qquad
\bw=(1,3,0,0),
\qquad
\bbeta=(0,1)^T.
\]
We take
\[
\bv_0=(0,0,-1,1)^T,
\qquad
B=
\begin{pmatrix}
  -1&1\\
  1&0\\
  1&-3\\
  -1&2
\end{pmatrix}.
\]
Then the negative support data in \Cref{app:negative-support-example-1} give 
\[
N=\{\{3\}=I_{\bzero},\{1\},\{1,3\},\{1,3,4\}\}
\] 
and $K_N=\emptyset$.
Thus $\widehat m_N(\bt)=t_3$ and $m_N(\bs)=(B\bs)_3$.

A direct computation gives $\widehat Q_N=\langle t_1t_j,\ t_3t_j\mid j=1,2,3,4\rangle$.
Hence
\begin{equation}
\label{eq:example1-Qhat-colon}
\widehat Q_N:\widehat m_N(\bt)=\widehat Q_N:t_3=\langle t_1,t_2,t_3,t_4\rangle.
\end{equation}
On the intrinsic side,
\begin{equation}
\label{eq:example1-P-colon}
P_N:m_N(\bs)=\langle(B\bs)_1,\ (B\bs)_2,\ (B\bs)_3,\ (B\bs)_4\rangle.
\end{equation}
Since $\Phi(t_j)=(B\bs)_j$ for $j=1,\dots,4$, \eqref{eq:example1-Qhat-colon} and \eqref{eq:example1-P-colon} imply \eqref{eq:examples-comparison-equality}.

This example also shows that the sufficient condition in
\Cref{prop:colon-y-general-N} is not necessary.

Indeed, from
\Cref{app:negative-support-example-1} one obtains $P_{N,I_{\bzero}}=\langle(B\bs)_2,\ (B\bs)_1(B\bs)_4\rangle$, whereas
$P_N:m_N(\bs)$ is the ideal in \eqref{eq:example1-P-colon}.
Moreover,
\[
\widehat Q_N:(t_3t_j)=\langle 1\rangle\neq\widehat Q_N:t_3\qquad(j=1,2,3,4),
\]
so the sufficient condition in \Cref{prop:square-general-N} is not necessary either.
\end{example}

\begin{example}[{\cite[Example 4.2.7]{SST}}]
\label{ex:2.2}
Let
\[
A=
\begin{pmatrix}
  1&1&1&1&1&1\\
  0&1&1&0&-1&-1\\
  -1&-1&0&1&1&0
\end{pmatrix},
\ \
\bw=(0,-1,-10,-10^2,-10^3,-10^4),
\]
and let $\bbeta=(1,0,0)^T$.  We take
\[
\bv_0=(-1,1,0,0,0,1)^T,
\qquad
B=
\begin{pmatrix}
  1&0&0\\
  -2&1&0\\
  2&-2&1\\
  -1&2&-2\\
  0&-1&2\\
  0&0&-1
\end{pmatrix}.
\]
The negative support data in \Cref{app:negative-support-example-2} give an ordered negative support family $N$ with $I_{\bzero}=\{1\}$ and $K_N=\emptyset$.
Thus $\widehat m_N(\bt)=t_1$ and $m_N(\bs)=s_1$.

In this example, $\widehat M_N=\langle t_1,\dots,t_6\rangle$ and a direct computation gives $\widehat Q_N=\langle t_it_j\mid 1\le i\le j\le 6\rangle$.
Hence
\begin{equation}
\label{eq:example2-Qhat-colon}
\widehat Q_N:\widehat m_N(\bt)
=
\widehat Q_N:t_1
=
\langle t_1,\dots,t_6\rangle.
\end{equation}
On the intrinsic side, one obtains $P_N=\langle s_1^2,\ s_2^2,\ s_3^2,\ s_1s_2,\ s_1s_3,\ s_2s_3\rangle$, and therefore
\begin{equation}
\label{eq:example2-P-colon}
P_N:m_N(\bs)
=
P_N:s_1
=
\langle s_1,s_2,s_3\rangle.
\end{equation}
Since $\Phi(\langle t_1,\dots,t_6\rangle)=\langle s_1,s_2,s_3\rangle$, \eqref{eq:example2-Qhat-colon} and \eqref{eq:example2-P-colon} imply \eqref{eq:examples-comparison-equality}.

Finally, this example shows that the product-intersection condition in
\Cref{cor:product-intersection-general-N} is not necessary.
Indeed, since $\widehat M_N=\langle t_1,\dots,t_6\rangle$,
we have
\[
\langle A\bt\rangle\widehat M_N\subsetneq\langle A\bt\rangle=\langle A\bt\rangle\cap\widehat M_N.
\]
Nevertheless,
\[
(\langle A\bt\rangle\widehat M_N+\widehat P_N):t_1=(\langle A\bt\rangle\cap\widehat M_N+\widehat P_N):t_1=\langle t_1,\dots,t_6\rangle.
\]
Thus the comparison equality \eqref{eq:examples-comparison-equality} may hold even when the product-intersection condition fails.
\end{example}

\begin{question}
\label{ques:intrinsic-not-sufficient}
Does there exist an ordered negative support family $N$ such that
\begin{equation}
\label{eq:open-question-comparison}
\Phi(\widehat Q_N:\widehat m_N(\bt))\neq P_N:m_N(\bs)?
\end{equation}
Equivalently, does there exist a case where the intrinsic perturbation method does not realize the full coefficient space obtained by the ambient perturbation method?
This question remains open even for the family $N=\cN_{\bv_0}$.
\end{question}

\appendix
\crefalias{section}{appendix}
\crefalias{subsection}{appendix}

\section{Negative support data for the examples}
\label{app:negative-support-computations}

In this appendix we record the negative support data used in
\Cref{sec:examples}.
We give the lattice parametrizations, the resulting ordered negative support families, and their complements.

\subsection{\texorpdfstring{\Cref{ex:2.1}}{Example 7.1}}
\label{app:negative-support-example-1}

For the data in \Cref{ex:2.1}, we use the basis
\[
B=
\begin{pmatrix}
  -1&1\\
  1&0\\
  1&-3\\
  -1&2
\end{pmatrix}
\]
of $L_{\C}$ and the fake exponent $\bv_0=(0,0,-1,1)^T$.
For $\lambda=(x,y)^T$, we have
\begin{equation}
\label{eq:example1-lattice-param}
\bv_0+B\lambda=(-x+y,\ x,\ -1+x-3y,\ 1-x+2y)^T.
\end{equation}
The weight condition is $\bw\cdot B\lambda=2x+y\ge0$.
Then
\[
\cN_{\bv_0}
=
\{\{3\},\{1\},\{1,3\},\{1,3,4\}\}.
\]
Thus, in \Cref{ex:2.1}, we take
\[
N=\cN_{\bv_0}
=
\{\{3\}=I_{\bzero},\{1\},\{1,3\},\{1,3,4\}\},
\qquad
K_N=\emptyset.
\]
The full set of negative supports occurring from
\eqref{eq:example1-lattice-param} is
\[
\NS_{\bv_0}
=
\{
\{1\},\{2\},\{3\},
\{1,3\},\{1,4\},\{2,3\},\{2,4\},
\{1,2,4\},\{1,3,4\}
\}.
\]
Therefore
\[
N^c
=
\NS_{\bv_0}\setminus N
=
\{\{1,4\},\{1,2,4\},\{2,4\},\{2\},\{2,3\}\}.
\]

\subsection{\texorpdfstring{\Cref{ex:2.2}}{Example 7.2}}
\label{app:negative-support-example-2}

For the data in \Cref{ex:2.2}, we use the basis
\[
B=
\begin{pmatrix}
  1&0&0\\
  -2&1&0\\
  2&-2&1\\
  -1&2&-2\\
  0&-1&2\\
  0&0&-1
\end{pmatrix}
\]
of $L_{\C}$ and the fake exponent
\[
\bv_0=(-1,1,0,0,0,1)^T.
\]
For $\lambda=(x,y,z)^T$, we have
\begin{equation}
\label{eq:example2-lattice-param}
\bv_0+B\lambda
=
(x-1,\ -2x+y+1,\ 2x-2y+z,\ -x+2y-2z,\ -y+2z,\ -z+1)^T.
\end{equation}
The weight condition is $\bw\cdot B\lambda=82x+819y+8190z\ge0$.
The ordered negative support family used in \Cref{ex:2.2} is
\[
\begin{aligned}
N=\{&
\{1\}=I_{\bzero},\{2\},\{3\},\{4\},\{5\},\{6\},
\{2,6\},\{3,6\},\{4,6\},\\
&\{1,3,6\},\{2,3,6\},\{2,5,6\},
\{3,4,6\},\{3,5,6\},\\
&\{1,3,5,6\},\{2,3,5,6\},\{2,4,5,6\}\}.
\end{aligned}
\]
In this example, $K_N=\emptyset$.
The complement of $N$ in $\NS_{\bv_0}$ is
\[
\begin{aligned}
N^c=\{&
\{1,3\},\{1,4\},\{1,5\},\{2,4\},\{2,5\},\{3,5\},\\
&\{1,2,4\},\{1,2,5\},\{1,3,4\},\{1,3,5\},\{1,4,5\},\\
&\{1,4,6\},\{2,3,5\},\{2,4,5\},\{2,4,6\},\\
&\{1,2,3,5\},\{1,2,4,5\},\{1,2,4,6\},\\
&\{1,3,4,5\},\{1,3,4,6\}\}.
\end{aligned}
\]

\bibliographystyle{amsplain}
\bibliography{refs}

@article{GGZ,
  author  = {Gelfand, I. M. and Graev, M. I. and Zelevinsky, A. V.},
  title   = {Holonomic systems of equations and series of hypergeometric type},
  journal = {Soviet Math. Dokl.},
  volume  = {36},
  year    = {1988},
  pages   = {5--10}
}

@article{GZK2,
  author  = {Gel'fand, I. M. and Zelevinsky, A. V. and Kapranov, M. M.},
  title   = {Hypergeometric functions and toral manifolds},
  journal = {Funct. Anal. Appl.},
  volume  = {23},
  year    = {1989},
  pages   = {94--106},
  note    = {English translation of Funktsional. Anal. i Prilozhen. 23 (1989), 12--26}
}

@article{GKZ,
  author  = {Gel'fand, I. M. and Kapranov, M. M. and Zelevinsky, A. V.},
  title   = {Generalized Euler integrals and {$A$}-hypergeometric functions},
  journal = {Adv. Math.},
  volume  = {84},
  year    = {1990},
  pages   = {255--271}
}

@article{Hotta,
  author  = {Hotta, R.},
  title   = {Equivariant {$D$}-modules},
  journal = {arXiv:math/9805021},
  year    = {1998}
}

@article{Log2,
  author  = {Okuyama, G. and Saito, M.},
  title   = {Logarithmic {$A$}-hypergeometric series {II}},
  journal = {Beitr. Algebra Geom.},
  volume  = {64},
  number  = {4},
  year    = {2023},
  pages   = {1057--1086},
  doi     = {10.1007/s13366-022-00669-5}
}

@article{LogFree,
  author  = {Saito, M.},
  title   = {Logarithm-free {$A$}-hypergeometric series},
  journal = {Duke Math. J.},
  volume  = {115},
  year    = {2002},
  pages   = {53--73}
}

@article{Log,
  author  = {Saito, M.},
  title   = {Logarithmic {$A$}-hypergeometric series},
  journal = {Int. J. Math.},
  volume  = {31},
  number  = {13},
  year    = {2020},
  pages   = {2050110},
  doi     = {10.1142/S0129167X20501104}
}

@book{SST,
  author    = {Saito, M. and Sturmfels, B. and Takayama, N.},
  title     = {Gr{\"o}bner Deformations of Hypergeometric Differential Equations},
  series    = {Algorithms and Computation in Mathematics},
  volume    = {6},
  publisher = {Springer},
  address   = {New York},
  year      = {2000}
}

@article{SW,
  author  = {Schulze, M. and Walther, U.},
  title   = {Irregularity of hypergeometric systems via slopes along coordinate subspaces},
  journal = {Duke Math. J.},
  volume  = {142},
  year    = {2008},
  pages   = {465--509}
}
 
\end{document}